\documentclass[aos,MSNbibl,nameyear,seceqn,dvips]{arximspdf}
\usepackage{dcolumn}
\usepackage{graphicx}

\doi{10.1214/13-AOS1116} 
\volume{41}
\issue{3}
\pubyear{2013}
\firstpage{1485}
\lastpage{1515}

\makeatletter
\newcommand{\rrvert}{\vert}
\newcommand{\llvert}{\vert}
\newcolumntype{d}[1]{D{.}{.}{#1}}
\newcommand{\eqref}[1]{(\ref{#1})}
\newcommand{\ran}{\operatorname{ran}}

\newcommand{\R}{\mathbb{R}}
\newcommand{\E}{\mathbb{E}}
\newcommand{\PP}{\mathbb{P}}
\newcommand{\HH}{\mathbb{H}}

\newcommand{\al}{\alpha}
\newcommand{\eps}{\varepsilon}
\newcommand{\om}{\omega}

\newcommand{\weak}{\stackrel{w}{\longrightarrow}}

\newcommand{\Wal}{\widetilde{\alpha}}

\newcommand{\Bga}{\overline{\gamma}}

\newcommand{\Var}{\operatorname{Var}}
\newcommand{\Cov}{\operatorname{Cov}}

\newcommand{\vect}[1]{\mathbf{#1}}

\newtheorem{prop}{Proposition}[section]
\newproclaim{ass}[prop]{Assumption}
\newtheorem{cor}[prop]{Corollary}
\newproclaim{defi}[prop]{Definition}
\newtheorem{lem}[prop]{Lemma}
\newtheorem{theo}[prop]{Theorem}
\newproclaim{rem}[prop]{Remark}
\makeatother

\begin{document}
\begin{frontmatter}

\title{Nonparametric inference on L\'evy measures and copulas\thanksref{T1}}
\runtitle{Inference on L\'evy measures and copulas}

\begin{aug}
\author{\fnms{Axel} \snm{B\"ucher}\corref{}\ead[label=e1]{axel.buecher@ruhr-uni-bochum.de}}
\and
\author{\fnms{Mathias} \snm{Vetter}\ead[label=e2]{mathias.vetter@ruhr-uni-bochum.de}}
\runauthor{A. B\"ucher and M. Vetter}
\thankstext{T1}{Supported by the collaborative
research center ``Statistical modeling of nonlinear dynamic processes''
(SFB 823) of the German Research Foundation (DFG).}
\affiliation{Ruhr-Universit\"at Bochum}

\address{Fakult\"at f\"ur Mathematik\\
Ruhr-Universit\"at Bochum\\
Universit\"atsstra\ss e 150\\
44780 Bochum\\
Germany \\
\printead{e1}\\
\phantom{E-mail:\ }\printead*{e2}}
\end{aug}

\received{\smonth{2} \syear{2013}}

%
\begin{abstract}
In this paper nonparametric methods to assess the multivariate L\'evy
measure are introduced. Starting from high-frequency observations of a
L\'evy process $\vect X$, we construct estimators for its tail
integrals and the Pareto--L\'evy copula and prove weak convergence of
these estimators in certain function spaces. Given $n$ observations of
increments over intervals of length $\Delta_n$, the rate of convergence
is~$k_n^{-1/2}$ for $k_n = n \Delta_n$ which is natural concerning
inference on the L\'evy measure. Besides extensions to nonequidistant
sampling schemes analytic properties of the Pareto--L\'evy copula
which, to the best of our knowledge, have not been mentioned before in
the literature are provided as well. We conclude with a short
simulation study on the performance of our estimators and apply them to
real data.
\end{abstract}

%


\end{frontmatter}
%
\section{Introduction} \label{Intro}

The modeling and estimation of dependencies is attracting an increasing
attention over the last decades in various fields of science like
mathematical finance, actuarial science or hydrology, among others.

In discrete time models, one of the most popular approaches is the
concept of copulas which allows one to separate the effects of
dependence of a random vector from its univariate marginal behavior. In
the bivariate case, the copula of a continuous random vector $(X,Y)$ is
the unique function $C\dvtx [0,1]^2\rightarrow[0,1]$ for which the identity
\[
\PP[X\le x, Y\le y] = C \bigl(\PP[X\le x], \PP[Y\le y] \bigr)
\]
holds for all $(x,y) \in\R^2$.
This formula, known as Sklar's theorem, is usually interpreted in the
way that the copula completely characterizes the stochastic dependence
between $X$ and $Y$ and hence represents the primary object of interest
for investigating dependencies. For introductions to the concept of
copulas in the aforementioned fields of science see \citet
{mcnfreemb2005}, \citet{freevald1998}, \citet{genefavr2007} and
references therein. The books of \citet{joe1997} and \citet{nelsen2006}
provide compendiums on the mathematical background and on various
parametric models. The huge amount of applications gave rise to a great
demand for statistical methods, of which semi- and nonparametric
estimation in discrete time i.i.d. models has been investigated in
\citet
{genghoriv1995}, \citet{ferradweg2004} and \citet{segers2012}, among others.

On the other hand, a huge number of models in applied stochastics
relies on an underlying process which is defined in continuous time. A
basic tool in this framework is the class of (multidimensional) L\'evy
processes which provides a flexible way to model empirically observed
behavior and includes prime examples such as Brownian motion and the
(compound) Poisson process. Statistical methods in this context
(including the somewhat more involved one of It\^o semimartingales)
depend on the nature of the observation schemes which are usually
classified as high-frequency and low-frequency setups. In both areas
the literature on nonparametrics has grown considerably over the last
decade. General overviews on high-frequency statistics can be found in
\citet{jacoprot2012} and \citet{myklzhan2012}. To mention only a few
approaches in detail which are close to our focus on jump processes we
refer to \citet{jacod2007}, \citet{figueroa2009}, \citet{aitsjaco2009},
\citet{bolltodo2011} and \citet{todotauc2011}. Seminal papers in the low
frequency setting are, for instance, due to \citet{neumreis2009} and
recently to \citet{nickreis2012}.

Our aim in this work is to combine both strands of the literature and
to provide nonparametric methods to estimate the dependence structure
of a multivariate L\'evy process. For the sake of brevity we will
concentrate on the bivariate case solely, but extensions to the general
$d$-dimensional setting are straightforward to obtain as well. Thus,
let $\vect X=(X^{(1)},X^{(2)})$ be a two-dimensional L\'evy process
with L\'evy--It\^o decomposition
%
%
\begin{equation}
\label{levy} \quad\vect X_t = \vect a t + \vect{B}_t + \int
_0^t \int_{\| \vect{u} \|
\leq
1} \vect{u} (
\mu- \bar{\mu}) (ds,d\vect{u}) + \int_0^t \int
_{\|
\vect
{u} \| > 1} \vect{u} \mu(ds,d\vect{u}),
\end{equation}
where $\vect a \in\mathbb R^2$ is a drift vector, $\vect B$ is a
bivariate Brownian motion with some covariance matrix $\Sigma$ and
$\mu
$ and $\bar{\mu}$ are the jump measure of the L\'evy process and its
compensator, respectively. It is well-known that the compensator takes
the form $\bar{\mu}(ds,d\vect{u})= ds\,  \nu(d\vect{u})$, where $\nu
$ is
the so-called L\'evy measure of $\vect X$. Given the choice of the
truncation function $h(\vect{u})= 1_{\{\| \vect{u} \| > 1\}}$, the law
of $\vect X$ is uniquely determined by the L\'evy triplet $(\vect a,
\Sigma, \nu)$.

As noted above, in the framework of statistics for stochastic processes
it is inevitable to lose some words on the underlying observation
scheme. We decide to work in a high-frequency setting which means in
the simplest case that at stage~$n$ one is able to observe one
realization of the process $\vect X$ at the equidistant times $i \Delta
_n$, $i=0, \ldots, n$, for a mesh $\Delta_n \to0$. A more general
setup including irregularly spaced data and asynchronous observations
will be provided as well.
Within the class of high-frequency settings a further distinction
regards the nature of the covered time horizon. Usually, we have either
$n\Delta_n = T$, corresponding to a finite time horizon (a trading day,
say), whereas $n \Delta_n \to\infty$ means that the process is
eventually observed on the entire time span $[0, \infty)$.

Due to the independence of the continuous part and the jump part of a
L\'evy process, the analysis of the stochastic nature of $\vect X$
canonically splits into inference on the covariance matrix $\Sigma$ and
inference on the jump measure~$\nu$.
However, estimation of the characteristics of the Brownian part of
$\vect X$ with or without additional jumps is well understood in the
high-frequency setup [among others, see \citet{jacod2008} for a thorough
theory on the behavior of more general It\^o semimartingales], so our
focus in this paper will be on the jump dependence of the two
components. In analogy to standard copulas for random vectors we will
employ a concept of a L\'evy copula to capture the dependence structure
within $\nu$ which dates back to \citet{conttank2004} and \citet
{kalltank2006}. We will follow a slightly different approach due to
\citet{resnklue2008} and \citet{ederklue2012}, however, and focus on
nonparametric methods to assess the closely related Pareto--L\'evy
copula. Consult also \citet{bolletal2012} for related work on jump
dependence using extreme value theory.

Besides parametric approaches to infer the (Pareto) L\'evy copula such
as \citet{esmaklue2011}, nonparametric methods in this area are hardly
available. To the best of our knowledge, the only concept is due to the
unpublished work of \citet{laeven2011} who constructs an estimator for
the L\'evy copula based on a limit representation involving ordinary copulas.
Some asymptotic properties are provided, but no explicit proof is
available. On the other hand, since the (Pareto) L\'evy copula captures
the tendency of the process to have joint (largely negative) jumps, the
need for reliable nonparametric estimators is evident from practice,
particularly with a view on finance. See, for instance, \citet
{boecklue2009} who model operational risk via Pareto--L\'evy copulas.
This convinces us that there is a clear gap in the literature which we
aim to fill in this work.

In contrast to Laeven's method, our approach will be based directly on
the defining relation of the Pareto--L\'evy copula $\Gamma$ which
involves tail integrals of both the L\'evy measure and its marginals.
For simplicity, we will focus on the spectrally positive case only,
that is, we assume that $\vect X$ has only positive jumps in both
directions. Equivalently, the L\'evy measure $\nu$ has support on
$[0,\infty)^2 \setminus\{(0,0)\}$, and $\Gamma$ will then naturally
be a function on the same space. In the case where all marginal tail
integrals have full range $[0,\infty]$, we obtain a representation of
$\Gamma$ as a functional of those, and we propose to estimate $\Gamma$
by using appropriate estimators for the tail integrals. It turns out
that in order to do so, we are forced to work in the high-frequency
setting with infinite time horizon, that is, $n \Delta_n \to\infty$.
Under some rather mild assumptions we are then able to prove weak
convergence of a suitably standardized version of $\hat\Gamma- \Gamma
$ in a certain function space, which will be our main result. As a
by-product, we obtain a Donsker theorem for the bivariate L\'evy
measure as well; a result which is similar in spirit to the recent work
of \citet{nickreis2012}, but in a high-frequency setting rather than a
low-frequency world.

The paper is organized as follows: Section~\ref{sec:plc} is devoted to
a brief discussion on jump dependence of bivariate L\'evy processes. We
summarize the concept of Pareto--L\'evy copulas and derive some of
their analytical properties. In Section~\ref{sec:est} we define
estimators for bivariate tail integrals, as well as for their
associated Pareto--L\'evy copulas. Weak convergence of these estimators
is discussed in Section~\ref{sec:weak}, while Section~\ref{sec:irr} is
devoted to deviations from the ideal sampling scheme. A~brief
discussion of our results, a small simulation study and a real data
example are provided in Section~\ref{sec:sim}, whereas some conclusions
are given in Section~\ref{sec:conc}. Finally, the main steps of the
proofs are postponed to Section~\ref{sec:aux}, while more technical
details are treated in a supplementary Appendix.

\section{Jump dependence and the Pareto--L\'evy copula} \label{sec:plc}
Suppose that we are given a bivariate L\'evy process $\vect X$ of the
form (\ref{levy}) where $\nu$ denotes its L\'evy measure. As already
stated in the \hyperref[Intro]{Introduction}, one assumption will be that $\nu$ has
support on $[0,\infty)^2 \setminus\{(0,0)\}$, which means that both
components of $\vect X$ only have positive jumps. This condition is for
notational convenience in first place, as we will see later that one
can follow a similar approach in order to estimate the jump dependence
in the other three quadrants as well.

Let us review some recent concepts of jump dependence. The basic
quantity in this framework is the bivariate tail integral $U$
associated with $\nu$, which, for the moment, will be defined as a
function from $[0,\infty]^2 \setminus\{(0,0)\}$ to $\R$ by
\[
U(\vect{x}) = \nu\bigl([x_1, \infty] \times[x_2,\infty]
\bigr), \qquad\vect{x}=(x_1,x_2).
\]
From the theory of L\'evy processes it is well known that this quantity
gives the average amount of jumps of $\vect X$ which fall into the
interval $[x_1, \infty] \times[x_2,\infty]$ during a time period of
length one. Since $\vect X$ has c\`adl\`ag paths, $U(\vect{x})$ is
necessarily finite. In the same way, we are able to introduce marginal
tail integrals. Precisely, let $U_i\dvtx [0,\infty] \to[0,\infty]$,
$i=1,2$, be defined via
\[
U_1(x_1) = \nu\bigl([x_1, \infty] \times
\mathbb R\bigr) \quad\mbox{and}\quad U_2(x_2) = \nu\bigl(\mathbb R
\times[x_2, \infty]\bigr).
\]
Again, $U_i(x_i)$ is finite for $x_i > 0$, but in the infinite activity
case we may have $U_i(0)=\infty$, $i=1,2$.

It is obvious that the entire information about $\nu$ is contained in
the tail integral~$U$. Therefore, just as for regular copulas, one
might be interested in splitting $U$ into several functions which are
related to the jump behavior of $\vect X$ in the marginals (naturally
given by the univariate tail integrals $U_i$) and a L\'evy copula $C$
which captures the specific tendency of $\vect X$ to have joint jumps.
Having this intuition in mind, \citet{conttank2004} provided the
following definition.

\begin{defi}
A bivariate L\'evy copula for L\'evy processes with positive jumps is a
function $C\dvtx [0, \infty]^2 \setminus\{(\infty,\infty)\} \to[0,
\infty
)$ which:
\begin{longlist}[(iii)]
\item[(i)] is grounded, that is, $C(x,0)=C(0,x)=0$ for all $x \in[0,
\infty]$;
\item[(ii)] has uniform margins, so $C(x,\infty)=C(\infty,x)=x$ for all
$x \in[0, \infty)$;
%
\item[(iii)] is 2-increasing, that is,
$C(x_1,x_2)-C(x_1,y_2)-C(y_1,x_2)+C(y_1,y_2) \geq0$ for all $x_1 \leq
y_1$ and $x_2 \leq y_2$.
\end{longlist}
\end{defi}

The main result on L\'evy copulas is a version of Sklar's famous
theorem which states that for each tail integral $U$ with marginals
$U_1$ and $U_2$ there exists a L\'evy copula $C$ such that
\[
U(\vect{x}) = C\bigl(U_1(x_1),U_2(x_2)
\bigr), \qquad\vect{x}=(x_1,x_2) \in [0,\infty]^2
\setminus\bigl\{(0,0)\bigr\},
\]
holds. Similarly to the usual copula, $C$ is uniquely defined on $\ran
U_1\times\ran U_2$ and, in particular, $C$ is globally unique if $\ran
U_i=[0,\infty]$ for $i=1,2$. This becomes our second main assumption.
In that case, we obtain a representation of $C$ via
\[
C(\vect{u}) = U\bigl(U_1^{-}(u_1),U_2^{-}(u_2)
\bigr),\qquad \vect{u}=(u_1,u_2) \in[0,\infty]^2
\setminus\bigl\{(\infty,\infty)\bigr\},
\]
where, for any $f\dvtx (0,\infty] \to[0,\infty)$ which is nonincreasing,
left-continuous and satisfies $f(\infty)=0$, $f^{-}\dvtx (0,\infty] \to
[0,\infty)$ denotes the generalized inverse function
%
%
\begin{equation}
\label{geni} f^{-}(z) = \inf\bigl\{x > 0 \dvtx f(x) \leq z\bigr\}
\end{equation}
and where we use the convention $f^-(0)=\infty$.

The inverse statement of Sklar's theorem is true as well: knowledge of
the marginals $U_i$ and the L\'evy copula $C$ completely determines $U$
and thus in turn~$\nu$. A drawback of the approach of \citet
{conttank2004} is, however, that $C$ is not a tail integral itself.
This is in contrast to the regular copula of a random vector which
couples marginal distribution functions and is a bivariate distribution
function on its own. This circumstance makes the interpretation of a L\'
evy copula quite difficult, and
for that reason it appears to be natural to focus on an alternative
notion of copula in this setting.

\begin{defi} \label{def:plc}
A bivariate Pareto--L\'evy copula for L\'evy processes with positive
jumps is a function $\Gamma\dvtx  [0, \infty]^2 \setminus\{(0,0)\} \to[0,
\infty)$ which:
\begin{longlist}[(iii)]
\item[(i)] is grounded, that is, $\Gamma(u,\infty)=\Gamma(\infty,u)=0$
for all $u \in(0, \infty]$;
\item[(ii)] has Pareto margins, so $\Gamma(u,0)=\Gamma(0,u)=1/u$ for
all $u \in(0, \infty]$;
\item[(iii)] is 2-increasing.
\end{longlist}
\end{defi}

As usual, we set $1/\infty= 0$ and vice versa.
Following \citet{ederklue2012}, Sklar's theorem now reads as follows:
given $U$ and its marginals, we have
\[
U(\vect{x}) = \Gamma \bigl(1/{U_1(x_1)},1/{U_2(x_2)}
\bigr),\qquad \vect {x}=(x_1,x_2) \in[0,\infty]^2
\setminus\bigl\{(0,0)\bigr\}
\]
for some unique Pareto--L\'evy copula $\Gamma$, and we obtain the relation
%
%
\begin{equation}
\label{PLC2}\qquad \Gamma(\vect{u}) = U \bigl(U_1^{-}(1/u_1),U_2^{-}(1/u_2)
\bigr), \qquad\vect{u}=(u_1,u_2) \in[0,\infty]^2
\setminus\bigl\{(0,0)\bigr\}.
\end{equation}
The difference to the approach of \citet{conttank2004} is that the
marginals of~$\Gamma$ correspond to Pareto tails, which are the tail
integrals of a 1-stable L\'evy process on the positive half line. Since
$\Gamma$ is 2-increasing as well, it is a simple task to deduce that it
satisfies the properties of a tail integral of a spectrally positive L\'
evy process as claimed.
Thus the Pareto--L\'evy copula allows for the interpretation that the
marginals of $\nu$ are standardized to the L\'evy measures of a
1-stable L\'evy process, which is similar in spirit to the ordinary
copula concept where marginals are standardized to uniform distributions.

Finally, we collect some basic properties of Pareto--L\'evy copulas,
some of which already have been stated in \citet{conttank2004} and
\citet
{kalltank2006} in the context of L\'evy copulas.

\begin{prop} \label{prop:plc}
Every Pareto L\'{e}vy copula $\Gamma$ has the following properties:
\begin{longlist}[(iii)]
\item[(i)]\textup{(Lipschitz continuity)}. $\llvert  \Gamma(\vect u) -
\Gamma
(\vect v) \rrvert  \le\llvert  \frac{1}{u_1} - \frac{1}{v_1} \rrvert  +
\llvert  \frac{1}{u_2} - \frac{1}{v_2} \rrvert $.\vspace*{1pt}
\item[(ii)]\textup{(Monotonicity)}. $\Gamma$ is 2-increasing and the
functions $\Gamma(u,\cdot)$ and $\Gamma(\cdot,u)$ are nonincreasing
for each fixed $u\ge0$.
\item[(iii)]\textup{(Fr\'echet-Hoeffding bounds}). $\Gamma_\perp\le\Gamma
\le\Gamma_\Vert $, where $\Gamma_\Vert (\vect u)= (u_1\vee
u_2)^{-1}$ and $\Gamma_\perp(\vect u)=u_1^{-1} 1_{\{ u_2=0\}} +
u_2^{-1} 1_{\{ u_1=0\}}$ denote the Pareto L\'{e}vy copulas
corresponding to perfect positive dependence and to independence, respectively.
\item[(iv)]\textup{(Partial derivatives)}. $\dot\Gamma_1(u_1,0)=-u_1^{-2}$
and $\dot\Gamma_1(u_1,\infty)=0$. For fixed $u_2\in(0,\infty)$, the
partial derivative $\dot\Gamma_1(u_1,u_2)$ exists for almost all
$u_1\in(0,\infty)$, and for such $u_1$ and $u_2$ we have
\[
0 \ge\dot\Gamma_1(u_1,u_2)
\ge-u_1^{-2}.
\]
Furthermore, the mapping $u_2\mapsto\dot\Gamma_1(u_1,u_2)$ is defined
and nondecreasing almost everywhere. Analogous results hold for the
partial derivative with respect to $u_2$.
\end{longlist}
\end{prop}

\section{Estimation of bivariate tail integrals and Pareto--L\'evy
copulas} \label{sec:est}

In the following we are interested in the construction of an estimator
$\hat\Gamma_n$ for $\Gamma$ which is based on relation (\ref{PLC2})
and empirical versions of the tail integrals $U$, $U_1$ and $U_2$. Such
estimators have, for instance, been discussed in \citet{figueroa2008} in
the univariate setting, and we will transfer them naturally to the
bivariate case.

Before we introduce these empirical versions, it turns out to be
convenient to slightly change the domain of $U$. Since by assumption no
negative jumps are involved, we have $
\nu([x_1, \infty] \times[0,\infty]) = \nu([x_1, \infty] \times
[-\infty,\infty])
$
for each $x_1>0$, and similarly for the second component. Therefore it
is equally well possible to define $U$ in the same way as before, but
as a function $U\dvtx \HH\to\R$, where
\[
\HH= (0,\infty]^2 \cup \bigl( \{-
\infty\} \times(0,\infty] \bigr) \cup \bigl( (0,\infty] \times\{ -\infty\} \bigr).
\]
Note that, on the stripes through $-\infty$, $U$ corresponds to the
marginal tail integrals $U_1$ and $U_2$, respectively.

Our estimator for the function $U$ will be defined on $\HH$ as well.
We set
%
%
\begin{equation}
\label{Un} 
U_n(\vect{x})=\frac{1}{k_n} \sum
_{j=1}^{n} 1_{\{\Delta_j^n X^{(1)}
\geq
x_1, \Delta_j^n X^{(2)} \geq x_2\}},\qquad \vect{x}=(x_1,x_2),
\end{equation}
where $k_n = n \Delta_n$ and $\Delta_j^n X^{(i)} = X^{(i)}_{j \Delta_n}
- X^{(i)}_{(j-1) \Delta_n}$ denotes the $j$th increment of $X^{(i)}$,
$i=1,2$. Having the role of the stripes through $-\infty$ in mind, we
obtain empirical versions of the univariate tail integrals through
%
%
\begin{equation}
\label{Uni} 
U_{n,1}(x_1)=U_n(x_1,-
\infty)=\frac{1}{k_n} \sum_{j=1}^{n}
1_{\{
\Delta
_j^n X^{(1)} \geq x_1\}},\qquad x_1\in(0,\infty],
\end{equation}
%
and analogously for $U_{n,2}$. Weak convergence of $U_n$ in an
appropriate function space is established in Theorem~\ref{prop1} below.

The underlying idea behind $U_n$ is rather natural, given the
interpretation of $U$ as the average number of jumps of a certain size
during the unit interval. Stationarity and indepedence of increments of
a L\'evy process ensure that the same behavior is to be expected over
intervals of arbitrary size, as long as $U$ is standardized
accordingly. Therefore, a canonical idea is to count joint large
increments of $X^{(1)}$ and $X^{(2)}$, as they indicate joint large
jumps over the corresponding time interval. This is precisely what
$U_n$ does.

\begin{rem} \label{rem:obs}  Several comments regarding the
conditions on the underlying sampling scheme are in order.
\begin{longlist}[(iii)]
\item[(i)] In order for $U_n$ to be consistent, it is necessary to be
in the high-frequency setting with infinite time horizon, that is, $k_n
\to\infty$. If we restrict ourselves to observations from a fixed time
interval $[0,T]$,\vadjust{\goodbreak} there are only finitely many jumps larger than a
given size, which is clearly not sufficient to draw inference on the
entire distribution of the jumps.
\item[(ii)] In reality, one usually neither sees both components of
$\vect X$ at the same time nor has a regular sampling scheme of
observations with distance $\Delta_n$. Therefore, the methods proposed
in this section need to be sharpened when dealing with a more general
setup including irregularly spaced data and asynchronous observations
which is done in Section~\ref{sec:irr} below. The same comment applies
concerning microstructure noise issues, in which case the distribution
of the increments depends heavily on the distribution of the noise
variables added.
\item[(iii)] The method of counting large increments in order to
identify jumps is a feature of observation schemes at high frequency.
In a low frequency world, with $\Delta_n = \Delta$ being fixed, one
rather observes a superposition of jumps and increments of the
continuous part and needs deconvolution techniques to distinguish
between both effects.
\end{longlist}
\end{rem}

In order to construct an empirical version of (\ref{PLC2}) recall the
definition of a generalized inverse function in \eqref{geni}.


\begin{defi}
Let $U$ be the tail integral of a bivariate L\'evy process with
positive jumps and $U_1$, $U_2$ be its marginal tail integrals.
Using their empirical versions (\ref{Un}) and (\ref{Uni}) we define,
for any $\vect{u}=(u_1,u_2) \in[0,\infty]^2\setminus\{ (0,0) \}$,
the \textit{empirical Pareto--L\'evy copula} as
%
%
\begin{equation}
\label{hatgam} 
\hat\Gamma_n(\vect{u}) =
U_n \bigl(\overline{U_{n,1}^-(1/u_1)},
\overline{U_{n,2}^-(1/u_2)} \bigr), 
\end{equation}
where $U_{n,i}^-$ is the generalized inverse function of $U_{n,i}$ as
defined in (\ref{geni}),
with the convention that $U_{n,i}^-(1/\infty)=U_{n,i}^-(0)=\infty$
and where $\overline{a}=a1_{\{a>0\}}-\infty1_{\{a=0\}}$ for some
$a\in
[0,\infty]$.
Finally, we set $U_n(-\infty,-\infty)=n/k_n$.
\end{defi}

\begin{rem}
In order to understand why $\overline{a}$ is introduced, suppose that
we are interested in estimating $\Gamma(u_1,0)$ even though it is known
to take the value $1/u_1$. Our estimator becomes
$U_n({U_{n,1}^-(1/u_1)}, -\infty)$ then, which is in general close to
$1/u_1$ due to the definition of $U_{n,1}$. On the other hand, if we
forget about $\overline{a}$, we obtain $U_n({U_{n,1}^-(1/u_1)}, 0)$
which only counts those increments of $\vect X$ where the first
component exceeds ${U_{n,1}^-(1/u_1)}$ and the second one is
nonnegative. Due to the existence of a Brownian part in $\vect X$,
however, we cannot expect these two estimators to be close, since a
nonnegligible number of increments in the second component is indeed
negative, and thus this estimator is considerably smaller than $\hat
\Gamma(u_1,0)$. 
\end{rem}

\begin{rem}
In the general case of arbitrary jumps, a similar construction allows
for the estimation of $\Gamma$ in the interior of each of the four
quadrants separately. Indeed, \citet{ederklue2012} give a general notion
of tail integrals and Pareto--L\'evy\vadjust{\goodbreak} copulas in their Definition 4, and
from Sklar's theorem in this context (which is their Theorem 1) we know
that the same relation as (\ref{PLC2}) holds for $\vect{u} \in
(\mathbb
R \setminus\{0\})^2$ and determines $\Gamma$ uniquely. For the sake
of brevity we dispense with the entire theory in this setting. 
\end{rem}

\section{Results on weak convergence} \label{sec:weak}

Our aim in this section is to prove results on weak convergence of both
estimators, and we begin with such a claim for $U_n$, as this theorem
is used later to show weak convergence of $\hat\Gamma_n$. Before we
come to the main theorems, let us briefly resume our assumptions on
$\nu
$ which mostly have already been given in the previous paragraphs.

\begin{ass} \label{ass}
Let $\vect X$ be a bivariate L\'evy process with the representation~(\ref{levy}). The following assumptions on $\nu$ are in order:
\begin{longlist}[(iii)]
\item[(i)] $\nu$ has support $[0, \infty)^2\setminus\{(0,0)\}$.
\item[(ii)] On this set it takes the form $\nu(d\vect u) = s(\vect u)
\,d\vect u$ for a positive L\'evy density $s$ which satisfies
\[
\sup_{u \in M_\eta} \bigl(\bigl|s(\vect u)\bigr| + \bigl\|\nabla s(\vect u)\bigr\|\bigr) <
\infty
\]
for any $\eta\in(0,\infty)^2$, where
\[
M_\eta= (\eta,\infty)^2 \cup\bigl(\{0\} \times(\eta,
\infty)\bigr) \cup \bigl((\eta,\infty) \times\{0\}\bigr),
\]
and $\nabla s$ denotes the gradient of $s$ on $(\eta,\infty)^2$ and the
univariate derivative on the stripes through 0, respectively.
\item[(iii)] $\nu$ has infinite activity, that is, $\nu([0, \infty)
\times[0,\infty)) = \infty$.
\end{longlist}
\end{ass}

Assumption~\ref{ass}(ii) had not been stated previously. It is used to
prove a second order condition regarding the difference between $U$ and
the expectation of $U_n$, for which we generalize a result due to \citet
{figuhoud2009} from the univariate setting to the multidimensional
case. Continuity and (strict) monotonicity of the marginal tail
integrals as claimed before are obvious consequences of it.

We begin with a result on weak convergence of $U_n$, and to this end we
have to define the function space on which the asymptotics take place.
Let $\mathcal B_\infty(\HH)$ be the space of all functions $f\dvtx \HH
\to
\mathbb R$ which are bounded on any subset of $\HH$ that is bounded
away from the origin and from the points $(-\infty,0)$ and $(0,-\infty
)$. We consider the metric inducing the topology of uniform convergence
on those subsets, defined by
\[
d(f,g) = \sum_{k=1}^\infty2^{-k} \bigl(
\| f-g \|_{T_k} \wedge1 \bigr),
\]
where $T_k = [1/k,\infty]^2 \cup(\{-\infty\} \times[1/k,\infty])
\cup
([1/k,\infty] \times\{-\infty\})$
and $\| f \|_{T_k} = \sup_{\vect u \in T_k} |f(\vect u)|$. This space
is a complete metric space, and a sequence converges in $\mathcal
B_\infty(\HH)$, if and only if it converges uniformly on each $T_k$.\vadjust{\goodbreak}

\begin{theo} \label{prop1}
Assume that $\vect X$ is a L\'evy process satisfying \textup{(i)} and \textup{(ii)} of
Assumption~\ref{ass}. If the observation scheme meets the conditions
%
%
\begin{equation}
\label{eq:cond} \Delta_n \to0,\qquad k_n \to\infty,\qquad
\sqrt{k_n} \Delta_n \to0,
\end{equation}
then we have
\[
\gamma_n(\vect{x})=\sqrt{k_n} \bigl\{ U_n(
\vect{x}) - U(\vect{x}) \bigr\} \weak{\mathbb B}(\vect{x})
\]
in $(\mathcal B_\infty(\HH),d)$, where ${\mathbb B}$ is a tight,
centered Gaussian process with covariance
\[
\E \bigl[{\mathbb B}(\vect{x}) {\mathbb B}(\vect{y}) \bigr] = U(\vect{x} \vee
\vect{y}) = U(x_1 \vee y_1, x_1 \vee
y_2).
\]
The sample paths of ${\mathbb B}$ are uniformly continuous on each
$T_k$ with respect to the pseudo distance
\[
\rho(\vect{x}, \vect{y}) = \E \bigl[ \bigl( {\mathbb B}(\vect{x}) - {\mathbb B}(
\vect{y}) \bigr)^2 \bigr]^{1/2} 
=
\bigl\llvert U(\vect{x}) - U(\vect{y}) \bigr\rrvert ^{1/2}.
\]
%
\end{theo}

For the proof of Theorem~\ref{prop1} the following lemma is extremely
useful. Its univariate version is a special case of a more general
result in \citet{figuhoud2009}.

\begin{lem} \label{lem1}
Suppose that \textup{(i)} and \textup{(ii)} of Assumption~\ref{ass} hold and let $\delta
> 0$ be fixed. Then there exist constants $K = K(\delta)$ and $t_0 =
t_0(\delta)$ such that the uniform bound
\[
\bigl\llvert \PP\bigl(X_t^{(1)} \geq x_1,
X_t^{(2)} \geq x_2\bigr) - t \nu
\bigl([x_1,\infty ) \times[x_2,\infty)\bigr) \bigr\rrvert
< K t^2
\]
holds for all $\vect x=(x_1,x_2)\in[\delta,\infty]^2 \cup ( \{
-\infty\} \times[\delta,\infty]  ) \cup ( [\delta,\infty]
\times\{ -\infty\}  )$ and $0 < t < t_0$.
\end{lem}

\begin{rem} \label{rem:bias}
Lemma~\ref{lem1} is used in the proof of Theorem~\ref{prop1} to show
that the bias of $U_n$ is of order $\Delta_n$.
It is this result which is responsible for the condition $\sqrt{k_n}
\Delta_n \to0$ in (\ref{eq:cond}) as the latter secures that this bias
term is negligible in the asymptotics. If one has $\sqrt{k_n} \Delta_n
\to c> 0$, instead, together with a stronger condition on
differentiability of $s(\vect u)$, a bias term will appear in Theorem
\ref{prop1}. Precisely, generalizing a result from \citet{figuhoud2009}
again, we obtain
\[
\PP\bigl(X_t^{(1)} \geq x_1,
X_t^{(2)} \geq x_2\bigr) = t \nu
\bigl([x_1,\infty) \times[x_2,\infty)\bigr) +
t^2 d_2(\vect x)/2 + O\bigl(t^{3}\bigr),
\]
uniformly in the same sense as above, with
\begin{eqnarray*}
d_2(\vect x) &=& 2 \biggl( \int_{x_2}^\infty
s(x_1,w_2) \,dw_2 a_\eps^1
+ \int_{x_1}^\infty s(w_1,x_2)
\,dw_1 a_\eps^2 \biggr)
\\
&&{}- \biggl(\int_{x_2}^\infty s_1(x_1,w_2)
\,dw_2 \Sigma_{11} + \int_{x_1}^\infty
s_2(w_1,x_2) \,dw_1
\Sigma_{22} - 2 s(x_1,x_2) \Sigma_{12}
\biggr)
\\
&&{}- 2 \int\int_0^1 (1-\beta) \biggl( \int
_{x_2 - \beta u_2}^\infty s_1(x_1-\beta
u_1,w_2) \,dw_2 u_1^2
\\
&&\hspace*{86pt}{} + \int_{x_1 - \beta
u_1}^\infty s_2(w_1,x_2-
\beta u_2) \,dw_1 u_2^2
\\
&&\hspace*{101pt}{} -2 s(x_1-\beta u_1,x_2-\beta
u_2) u_1u_2 \biggr) \,d\beta
\nu_\eps(d\vect u)
\\
&&{}+ \int \int s(\vect z) \bar c_{\eps}(\vect z) s(\vect w) \bar
c_{\eps}(\vect w) 1_{\{\vect w+ \vect z \geq x\}} \,d \vect z\, d \vect w - 2 \lambda
_\eps \nu\bigl([\vect x,\infty)\bigr).
\end{eqnarray*}
Here, $\eps>0$ denotes an auxiliary variable used to simplify the
expression, and for all unexplained notation we refer to the proof of
Lemma~\ref{lem1}. Under the condition $\sqrt{k_n} \Delta_n \to c> 0$
above, it can be shown that the asymptotic bias in Theorem~\ref{prop1}
takes the form $cd_2(\vect x)/2$. Note that it is possible to derive a
representation for $d_2(\vect x)$ independently of $\eps$, which takes
an even more complicated form. This expression is given explicitly in
the supplementary material, alongside with a sketch of a proof.
\end{rem}

Before we come to the result on $\hat\Gamma_n$, let us introduce an
oracle estimator for $\Gamma$. We set, for any $\vect{u}=(u_1,u_2)
\in
[0,\infty]^2\setminus\{ (0,0) \}$,
\[
\widetilde\Gamma_n(\vect{u}) = U_n \bigl(U_{1}^{-}(1/u_1),U_{2}^{-}(1/u_2)
\bigr),
\]
which means that we replace the inverses of the empirical marginal tail
integrals by the unobservable true ones. Thanks to Theorem~\ref{prop1}
we obtain weak convergence of a restricted version of this intermediate
estimator in the space $\mathcal B_\infty((0,\infty]^2)$ of all real
functions on $(0,\infty]^2$ that are bounded on sets which are bounded
away from the origin. In a similar sprit as before, we equip this space
with the metric $d(f,g)=\sum_{k=1}^\infty2^{-k}  ( \| f-g \|_{T_k}
\wedge1  )$, where $T_k=[1/k,\infty]^2$. Setting $\vect
{x}=(U_1^{-}(1/u_1), U_2^{-}(1/u_2))$ and observing that $U_i^{-}(k)
\ge k'>0$, the continuous mapping theorem immediately yields the
following result.

\begin{cor}
Under the conditions of Theorem~\ref{prop1} we have
\[
\Wal_n(\vect{u})=\sqrt{k_n} \bigl( \widetilde
\Gamma_n(\vect{u}) - \Gamma(\vect{u}) \bigr) \weak\mathbb B
\bigl(U_1^{-}(1/u_1), U_2^{-}(1/u_2)
\bigr) 
\]
in $(\mathcal B_\infty((0,\infty]^2),d)$, with $\mathbb B$ as defined
in Theorem~\ref{prop1}.
\end{cor}

From a statistical point of view, there is no loss in information when
estimating $\Gamma(\vect u)$ on $(0,\infty]^2$ instead of the entire
domain $[0,\infty]^2\setminus\{ (0,0) \}$, since a Pareto--L\'evy
copula is grounded by definition and thus known on stripes through 0.
This remark remains valid for the final result of this section as well,
which is on weak convergence of the estimator $\hat\Gamma_n(\vect{u})$.

\begin{theo} \label{theo1}
Assume that $\vect X$ is a L\'evy process satisfying Assumption~\ref
{ass}. If~(\ref{eq:cond}) holds, then we have
\[
\al_n(\vect{u})=\sqrt{k_n} \bigl( \hat
\Gamma_n(\vect{u}) - \Gamma (\vect{u}) \bigr) \weak\mathbb G(\vect{u})
\]
in $(\mathcal B_\infty((0,\infty]^2),d)$. Here, the process $\mathbb G$
is defined as
%
%
\begin{equation}
\label{eq:gp} \mathbb G(\vect u) = \widetilde{\mathbb G}(\vect{u}) +
u_1^2 \dot \Gamma_1(\vect u) \widetilde{
\mathbb G}(u_1,-\infty) + u_2^2 \dot
\Gamma_2(\vect u) \widetilde{\mathbb G}(-\infty,u_2),
\end{equation}
where $\widetilde{\mathbb G}$ denotes a tight centered Gaussian field
on $\HH$ with covariance structure
\[
\E \bigl[\widetilde{\mathbb G}(\vect{u}) \widetilde{\mathbb G}(\vect {v}) \bigr]
= \Gamma(\vect{u} \vee\vect{v}) =\Gamma(u_1 \vee v_1,
u_2 \vee v_2)
\]
using the convention $\Gamma(u,-\infty) = \Gamma(-\infty,u) =1/u$. The
sample paths of $\widetilde{\mathbb G}$ are uniformly continuous on
each $T_k$ with respect to the pseudo distance
\[
\rho(\vect u,\vect v) = \E \bigl[ \bigl( \widetilde{\mathbb G}(\vect u) -
\widetilde{\mathbb G}(\vect v) \bigr)^2 \bigr]^{1/2}
=\bigl\llvert \Gamma(\vect u) - \Gamma(\vect v) \bigr
\rrvert ^{1/2}.
\]
\end{theo}

If both coordinates of $\vect u$ are distinct from $\infty$, then
$\dot
\Gamma_i(\vect u)$ exists as a consequence of (\ref{PLC2}) and
Assumption~\ref{ass}, and $\mathbb G(\vect u)$ is well-defined. On the
other hand, if one of the components equals $\infty$, we have
$\widetilde{\mathbb G}(\vect u)=0$ almost surely; and also $\dot
\Gamma
_1(u_1, \infty)=0$ and $\dot\Gamma_2(\infty,u_2)=0$ from Proposition
\ref{prop:plc}. Hence, the right-hand side of \eqref{eq:gp} is well
defined as well, and we have $\mathbb G(\vect u)=0$ almost surely in
this case.


%

\section{Deviations from the regular setting} \label{sec:irr}

Up to now, the results in this paper have been shown in the ideal
setting of observing a L\'evy process at equidistant and synchronous
times. With a view to applications, it is obvious that these
assumptions are not realistic in practice.
For example, due to the stylized facts of financial time series, pure
L\'evy processes possessing independent increments are too restrictive
for the modeling of financial time series. Also, multiple stock prices
are usually traded at different time points, contradicting our assumed
observation scheme.
Ways to overcome the latter problem have so far only been discussed in
the context of volatility estimation [a remarkable exception is \citet
{comtgeno2010}], and it is known that limit theorems usually differ
from the ones for regular observation times and are obtained under
quite restrictive assumptions only, particularly when the time points
are asynchronous. See, for example, \citeauthor{aitsmykl2003} (\citeyear{aitsmykl2003,aitsmykl2004}),
\citeauthor{hayayosh2005} (\citeyear{hayayosh2005,hayayosh2008}), \citet{myklzhan2009} or \citet{hayjacyos2011}.

Our aim in this section is to develop a concise theory which allows for
a consistent estimation of the L\'evy measure also in case of irregular
observations or when one observes a model with time-varying drift and
diffusion part. We will focus in particular on inference on the
distribution function $U(\vect x)$,\vadjust{\goodbreak} since weak convergence of the
empirical Pareto--L\'evy copula process carries over using the same
results on Hadamard differentiability as in the proof of Theorem \ref
{theo1}. Furthermore, we will indicate how results change in case of
microstructure noise and why standard methods for diffusions do not
carry over to our setting.

\subsection{Nonequidistant sampling schemes} \label{sec:irr1}
The first generalization regards the assumption of an equidistant
sampling scheme which we no longer assume to hold. Instead, suppose
from now on that the observation times are given by deterministic~$t_j^n$,
$j=1, \ldots, m_n$, where we set $k_n = t_{m_n}^n$ as before.
The results carry over to the case of random sampling as well, at least
if the observations times are independent of $\vect X$.

Still, the most natural estimator appears to be given by counting joint
large increments, which results in setting
%
%
\begin{equation}
\label{Vn} V_n(\vect{x})=\frac{1}{k_n} \sum
_{j=1}^{m_n} 1_{\{\Delta_j^n X^{(1)}
\geq x_1, \Delta_j^n X^{(2)} \geq x_2\}},
\end{equation}
where we have defined $\Delta_j^n X^{(i)} = X^{(i)}_{t_j^n} -
X^{(i)}_{t_{j-1}^n}$ (with $t_0^n=0$) in the same spirit as before.
Finally, let $\pi_n = \max_{j=1, \ldots, m_n} (t_j^n - t_{j-1}^n)$.
%
\begin{theo} \label{theo2}
Assume that $\vect X$ is a L\'evy process satisfying \textup{(i)} and \textup{(ii)} of
Assumption~\ref{ass}.
If furthermore the observation scheme meets the conditions
%
%
\begin{equation}
\label{samp1:cond} k_n \to\infty,\qquad \pi_n \to0,\qquad
\frac{1}{\sqrt{k_n}} \sum_{j=1}^{m_n}
\bigl(t_j^n - t_{j-1}^n
\bigr)^2 \to0,
\end{equation}
then we have $\sqrt{k_n}  \{ V_n(\vect{x}) - U(\vect{x})
\}
\weak{\mathbb B}(\vect{x})$ in $(\mathcal B_\infty(\HH),d)$, where
${\mathbb B}$ is the same Gaussian process as in Theorem~\ref{prop1}.
\end{theo}
%
\begin{rem} Even though the assumption $\pi_n \to0$ ensures that
we remain in a genuine high-frequency setting, it is not necessary in
general and only used here to simplify the proofs. Suppose, for
example, that we have $t_1^n = 1$ for all $n$, whence it is obviously
not possible to infer the L\'evy measure consistently over the interval
$[0,1]$. This, however, does not affect the validity of Theorem \ref
{theo2} per se, since a single increment of $\vect X$ contributes with
either $1/k_n$ or zero to $V_n$, so its influence is negligible in the
asymptotics, and as $k_n \to\infty$ holds, one can equally well use
the observations over $[1,k_n]$ to estimate the function $U$. 
\end{rem}
%
\subsection{Asynchronous sampling schemes} \label{sec:irr2}
Suppose now that both components of~$\vect X$ are observed at different
time stamps. We call $r_j^n$, $j=1, \ldots, m^1_n$, the series of
observations times connected with the process $X^{(1)}$, whereas
$s_\ell
^n$, $\ell=1, \ldots, m^2_n$, belongs to the second\vadjust{\goodbreak} component~$X^{(2)}$.
For simplicity only, we assume that the endpoints coincide,
that is $r_{m_n^1} = s_{m_n^2} = k_n$. As before, we denote with $\pi
_n^i$ the mesh of the $i$th time series.

In this situation it is less obvious how to count the number of joint
large jumps, as increments over $X^{(1)}$ and $X^{(2)}$ are in general
never computed over the same time intervals. In spirit of \citet
{hayayosh2005}, however, it appears reasonable to construct a na\"ive
estimator from counting those pairs of large increments which are
computed over at least overlapping intervals. Of course, in this case
it may happen that jumps at different times (but close nearby) are
treated as joint ones. For this reason it is important to assume
additional properties of the jump measures which make such an event unlikely.
%
\begin{ass} \label{ass2}
Let $\nu^i(dx)$ be the univariate L\'evy measures for $i=1,2$. We
assume that $\nu^i(dx) = \nu_1^i(dx) + \nu_2^i(dx)$ for mutually
singular measures $\nu_1^i$ and~$\nu_2^i$, given by
\[
\nu_1^i(dx) = a_i \frac{1+|x|^{\gamma_i}f_i(x)}{|x|^{1+\beta_i}}
1_{\{0
< x \leq1\}} \,dx
\]
with $0 \leq f_i(x) \leq K$, $\gamma_i \geq0$, $a_i \geq0$ and $\nu
_2^i$ such that $\nu_2^i(dx) = s^i(x) \,dx$ for some L\'evy density $s^i$
such that $\int(|x|^{\beta'_i} \wedge1) s^i(x) \,dx < \infty$, where $0
\leq\beta'_i < \beta_i < 2$.
\end{ass}

This condition means basically that the behavior of small jumps in both
components is similar to the one of $\beta_i$-stable processes. Such
assumptions are often used in high-frequency statistics for jump
processes; cf., for example, \citet{aitsjaco2009} and related work. Our
estimator for $U(\vect x)$ now reads as follows:
%
%
\begin{equation}
\label{Wn} W_n(\vect{x})=\frac{1}{k_n} \sum
_{j=1}^{m_n^1}\sum_{\ell=1}^{m_n^2}
1_{\{\Delta_j^n X^{(1)} \geq x_1, \Delta_\ell^n X^{(2)} \geq x_2\}} 1_{\{(r_{j-1}^n,r_j^n] \cap(s_{\ell-1}^n,s_{\ell}^n] \neq\varnothing
\}}.
\end{equation}
Under slightly more restrictive conditions on the sampling scheme than
before, we obtain the following result on weak convergence.
%
\begin{theo} \label{theo3}
Let $\vect X$ be a L\'evy process satisfying Assumptions~\ref{ass}\textup{(i)}
and \textup{(ii)} and Assumption~\ref{ass2}. Setting $\beta= \max(\beta_1,
\beta
_2)$ assume further that the sampling scheme satisfies $\max(\pi_n^1,
\pi_n^2) \to0$ as well as
\[
\frac{1}{\sqrt{k_n}} \sum_{j=1}^{m_n^1}
\bigl(r_j^n - r_{j-1}^n
\bigr)^{
{(\beta
+2)}/{(\beta+1)}} \to0 \qquad\mbox{if } \beta> 1
\]
and
\[
\frac{1}{\sqrt{k_n}} \sum_{j=1}^{m_n^1}
\bigl(r_j^n - r_{j-1}^n
\bigr)^{3/2-\delta
} \to0 \qquad\mbox{if } \beta\leq1
\]
for some $\delta\in(0,1/2)$, and similarly for the increments
involving $s_{\ell}^n$. Then we have $\sqrt{k_n}  \{ W_n(\vect{x})
- U(\vect{x})  \} \weak{\mathbb B}(\vect{x})$ in $(\mathcal
B_\infty(\HH),d)$, where ${\mathbb B}$ is the same Gaussian process as
in Theorem~\ref{prop1}.
\end{theo}

\begin{rem}
It is remarkable that even though methods similar to the case of
volatility estimation work for inference on $U$ as well, the results
look quite different in this setting: in case of irregular
observations, only very few additional assumptions on the observation
scheme are necessary, which is in contrast to the restrictive
conditions of \citet{hayayosh2008} regarding covolatility. This happens,
however, at the cost of additional assumptions on the structure of the
underlying L\'evy process. When dealing with irregular sampling times,
we obtain the same central limit theorem as for equidistant ones. This
has again no direct connection to volatility estimation, as the
corresponding result in \citet{myklzhan2009} comes with a different variance.
\end{rem}

\begin{rem}
There is a variety of models which satisfy Assumptions~\ref{ass} and
\ref{ass2}. Among the simplest are stable processes in both components,
coupled by some L\'evy copula which is twice continuously
differentiable away from the origin (e.g., the Clayton one). 
\end{rem}

\subsection{Observing semimartingales} \label{Semi}
In this section we discuss briefly deviations from assumption (\ref
{levy}) on the observed process. Suppose that the underlying process is
an It\^o semimartingale with the representation
\begin{eqnarray*}
\vect X_t &=& \int_0^t \vect
a_s \,ds + \int_0^t
\sigma_s \,d\vect B_s
\\
&&{}+ \int_0^t \int_{\| \vect{u} \| \leq1}
\vect{u} (\mu- \bar{\mu }) (ds,d\vect{u}) + \int_0^t
\int_{\| \vect{u} \| > 1} \vect{u} \mu (ds,d\vect{u})
\end{eqnarray*}
instead, where $\vect a \in\R^2$ and $\sigma\in\R^{2 \times2}$ are
bounded and left-continuous processes. Recall that all central limit
theorems proposed before deal with sums of large increments of $\vect
X$, interpreted as coming from large jumps over the same time
intervals. This intuition is based on the fact that the probability of
the continuous L\'evy part to become large over small time intervals is
exponentially small, and therefore it is likely that these claims
remain valid under weaker conditions on $\vect X$ as well. The
following theorem states that this is indeed the case, if we assume
that all sampling schemes satisfy similar growth conditions as in
Theorem~\ref{theo3}.
%
\begin{theo} \label{theo4}
Let $\vect X$ be an It\^o semimartingale as above, and assume that the
respective assumptions on the sampling schemes from
Theorems~\ref{theo1}, \ref{theo2} and~\ref{theo3} are satisfied. Assume further that
the sampling schemes satisfy
\[
\sqrt{k_n} \Delta_n^{1/2-\delta} \to0 \quad\mbox{or}\quad
\frac
{1}{\sqrt{k_n}} \sum_{j=1}^{m_n}
\bigl(t_j^n - t_{j-1}^n
\bigr)^{3/2-\delta} \to0
\]
for Theorems~\ref{theo1} or~\ref{theo2}, respectively, and for
some $\delta\in(0,1/2)$. Then, the weak convergence results of the
respective theorems hold.
\end{theo}

The proof relies on replacing increments of $\vect X$ by corresponding
increments of a pure jump L\'evy process in order to apply theorems on
weak convergence based on sums of independent observations. At this
stage, the additional conditions on the sampling scheme come into play.
Weaker conditions might be sufficient here, if one was able to prove
weak convergence based on some type of conditional independence
instead, then using different approximations for increments of $\vect X$.

\subsection{Microstructure noise} \label{Noise}
An important issue regarding high-frequency data is the presence of
microstructure noise, in which case one does not observe the plain L\'
evy process, but
\[
Z^{(j)}_{i\Delta_n} = X^{(j)}_{i\Delta_n} +
V^{(j)}_{i \Delta_n},
\]
where the $V^{(j)}_{i \Delta_n}$, $i=1, \ldots, n$, $j=1,2$, are
i.i.d. processes, independent of $X$, which satisfy $\E[V^{(j)}_{i
\Delta_n}] = 0$ and have moments of all order. In this case, our
estimator for the L\'evy measure does not work anymore. To see this,
let us for simplicity stick to the univariate setting. We have
\[
U_n(x) = \frac{1}{k_n} \sum_{i=1}^n
1_{\{\Delta_i^n Z \geq x\}}
\]
now, and if there is a positive probability of $P(\Delta_i^n V \geq
x)$, then $U_n(x)$ will behave like $\Delta_n^{-1} P(\Delta_i^n V
\geq
x)$ which diverges to infinity. On the other hand, if $P(\Delta_i^n V
\geq x) = 0$, then $U_n(x)$ is still bounded, but will rather estimate
a convolution of jumps and noise than the plain jump measure.

Microstructure noise issues are well understood in the context of
diffusions. However, it appears that none of the standard methods for
diffusion processes [let us mention the multiscale approach by \citet
{zhangetal2005} and the kernel-based one due to \citet{barndetal2008}]
can directly be applied in our context. Even the pre-averaging approach
by \citet{jacodetal2009}, which provides a general concept for
diminishing the influence of the noise by retaining information about
increments of $X$, fails when one is interested in the L\'evy measure.
In the following, we will briefly discuss these issues.

For an auxiliary sequence $l_n$ and some piecewise differentiable
function $g$ on $[0,1]$ with $g(0)=g(1)=0$,\vadjust{\goodbreak} \citet{jacodetal2009}
discuss $\widetilde{Z}_i^n = \widetilde{X}_i^n + \widetilde{V}_i^n$,
where, for an arbitrary $Y$ and $i=0, \ldots, n-l_n$,
\[
\widetilde{Y}_i^n = \sum_{j=1}^{l_n}
g(j/l_n) \Delta_{i+j}^n Y.
\]
While $\widetilde{X}_i^n$ can be seen as some kind of generalized
increment which still bears similar information as the plain $\Delta
_i^n X$, we have
\[
\widetilde{V}_i^n = -g(1/l_n)
V_{i \Delta_n} + \sum_{j=1}^{l_n-1}
\bigl(g(j/l_n)-g\bigl((j+1)/l_n\bigr)\bigr)
V_{(i+j) \Delta_n} = O_p \bigl(l_n^{-1/2} \bigr).
\]
Here we use both piecewise differentiability of $g$ and the assumptions
concerning $g$ on the boundary of $[0,1]$. Therefore the larger $l_n$
becomes, the less important is the contamination by noise.

For this reason, estimation based on pre-averaging usually works in the
way that one proceeds as usual, but replaces the standard estimators by
ones based on $\widetilde{Z}_i^n$. If $l_n$ is rather large compared
with $n$, it is reasonable to replace $\widetilde{Z}_i^n$ with
$\widetilde{X}_i^n$ in the asymptotics and thus to recover full
information of $X$. In our setting this would lead to an estimator of
the form
\[
\widetilde{U}_n(x) = \frac{1}{k_n} \sum
_{i=0}^{\lfloor n/l_n \rfloor
-1} 1_{\{\widetilde{Z}_{il_n}^n \geq x\}},
\]
where the $\widetilde{Z}_{il_n}^n$ are computed over nonoverlapping
intervals to retain i.i.d. terms in the sum. As noted above, for large
$l_n$ it is equally well possible to discuss $\widetilde{U}_n^X(x)$
which is defined similarly to $\widetilde{U}_n(x)$ above, but using
$\widetilde{X}_{il_n}^n$ instead of $\widetilde{Z}_{il_n}^n$. This
procedure, however, does not result in an estimator for the L\'evy
distribution function $U(y)$. The reason is that the leading term in an
expansion of $P(\widetilde{X}_{0}^n \geq y)$ is due to a single large
jump within $[0,l_n\Delta_n]$. In this case, its contribution to
$\widetilde{X}_{0}^n$ depends on its exact position within the
interval, as it has to be standardized by $g(j/l_n)$ accordingly. For
example, if the jump occurs in the small interval $[(j-1)\Delta_n,
j\Delta_n]$, it is not important whether the jump is larger than $y$,
but whether it is larger than $y/g(j/l_n)$. Since the jump time is
uniformly distributed, one can show formally that
\begin{eqnarray*}
P\bigl(\widetilde{X}_{il_n}^n \geq y\bigr) &=&
\Delta_n \sum_{j=1}^{l_n} U
\bigl(y/g(j/l_n)\bigr) + O\bigl(l_n^2
\Delta_n^2\bigr)
\\
&=& l_n \Delta_n \int_0^1
U\bigl(y/g(x)\bigr) \,dx + o(l_n \Delta_n),
\end{eqnarray*}
which proves that $\widetilde{U}_n^X(y)$ converges to $\int_0^1
U(y/g(x)) \,dx$ and is therefore not consistent for $U(y)$. Construction
of a consistent estimator for the L\'evy measure in case of
microstructure noise thus seems to be a challenging topic for future research.

\section{Discussion, simulations and an illustration} \label{sec:sim}

\subsection{An asymptotic comparison}

Suppose a statistician has knowledge of the marginal tail integrals. In
this case, the results in Section~\ref{sec:weak} provide two
competitive asymptotically unbiased estimators for the Pareto--L\'evy
copula, namely the oracle estimator $\widetilde\Gamma_n$ exploiting
knowledge of the marginals and the empirical Pareto--L\'evy copula
$\hat\Gamma_n$ ignoring this additional information. The following
proposition gives a partial answer to the question of which estimator
is (asymptotically) preferable. Perhaps surprisingly, ignoring the
additional knowledge decreases the asymptotic variance under certain
growth conditions on $\Gamma$. A similar observation has recently been
made in the context of copula estimation; see \citet{genesege2010}.

\begin{prop} \label{prop:comp}
Suppose that the Pareto--L\'evy copula $\Gamma$ has continuous
first-order partial derivatives and that the functions
%
%
\begin{equation}
\label{eq:mon} u_1\mapsto u_1 \Gamma(u_1,u_2)=
\frac{\Gamma(u_1,u_2)}{\Gamma(u_1,0)}, \qquad u_2\mapsto u_2 \Gamma(u_1,u_2)
= \frac{\Gamma(u_1,u_2)}{\Gamma(0,u_2)}
\end{equation}
are nondecreasing for fixed $u_2\in(0,\infty]$ and $u_1\in(0,\infty]$,
respectively. Then the Gaussian fields $\mathbb G$ and $\widetilde
{\mathbb G}$ satisfy the inequality
%
\[
\Cov\bigl\{ \mathbb G(\vect u), \mathbb G(\vect v) \bigr\} \le\Cov\bigl\{
\widetilde {\mathbb G}(\vect u), \widetilde{\mathbb G}(\vect v) \bigr\}
\]
for all $\vect u, \vect v \in(0,\infty]^2$. Particularly, $\Var\{
\mathbb G(\vect u) \} \le\Var\{ \widetilde{\mathbb G}(\vect u) \}$.
\end{prop}

\begin{figure}

\includegraphics{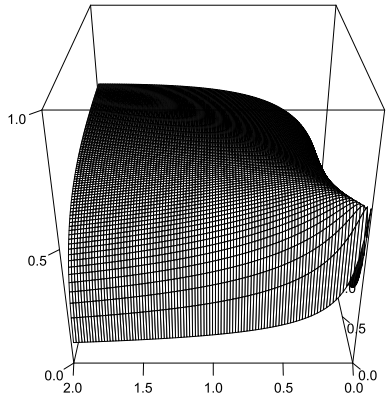}

\caption{The graph of the asymptotic relative efficiency of
$\widetilde
\Gamma_n$ to $\hat\Gamma_n$
for the Clayton Pareto--L\'evy copula with $\theta=0.5$.} \label{fig:levels}
\end{figure}

Under the assumptions of Proposition~\ref{prop:comp} the condition in
\eqref{eq:mon} is equivalent to
\[
u_1 \dot\Gamma_1(\vect u) + \Gamma(\vect u) \ge0,\qquad
u_2 \dot \Gamma_2(\vect u) + \Gamma(\vect u) \ge0
\]
for each $\vect u=(u_1,u_2)\in(0,\infty]^2$, which is easily accessible
for most parametric classes of Pareto--L\'evy copulas. For instance,
for the Clayton Pareto--L\'evy copula given by
\[
\Gamma(\vect u) = \bigl(u_1^\theta+ u_2^\theta
\bigr)^{-1/\theta}
\]
we have
\[
u_1 \dot\Gamma_1(\vect u) + \Gamma(\vect u) =
\bigl(u_1^\theta +u_2^\theta
\bigr)^{-1/\theta-1} u_2^\theta, \qquad u_2 \dot
\Gamma_2(\vect u) + \Gamma (\vect u) = \bigl(u_1^\theta+u_2^\theta
\bigr)^{-1/\theta-1} u_1^\theta,
\]
which is readily seen to be nonnegative. In Figure~\ref{fig:levels} we
depict the graph of the asymptotic relative efficiency
\[
[0,2]^2 \to[0,\infty), \qquad\vect u \mapsto{\Var\bigl\{ \mathbb G( \vect
u) \bigr\} }/{ \Var\bigl\{ \widetilde{ \mathbb G}(\vect u) \bigr\}}
\]
of the oracle estimator $\widetilde\Gamma_n$ to the empirical
Pareto--L\'evy copula $\hat\Gamma_n$ for $\vect u \in[0,2]^2$. The
Clayton parameter is chosen as $\theta=0.5$. Close to the axis the
relative efficiency decreases to $0$, while the maximal relative
efficiency is attained on the diagonal with a value of $21/32\approx
0.656$. Even in this best case, the difference is seen to be substantial.

\subsection{Simulation study}\label{subsec:sim}

In order to obtain an impression on the performance of the asymptotic
results stated in the previous section we will discuss some finite
sample properties concerning Theorems~\ref{prop1} and~\ref
{theo1}. In both cases, the setting is as follows: we simulate
(essentially) two $1/2$ stable subordinators, that is, both tail
integrals are given by $U_i(x) = (\pi x)^{-1/2}$, which are coupled by
a Clayton--Pareto--L\'evy copula with $\theta= 1/2$. Sometimes, we add
two independent Brownian motions with variance $1/2$ each, and
sometimes, we assume to observe the pure jump processes only.

Recall that the rate of convergence is $k_n^{-1/2}$ [which, in light of
the results in \citet{figuhoud2009}, appears to be natural in the
context of estimating the L\'evy measure]. Hence, a larger $k_n$
suggests a better approximation by the limiting Gaussian process,
whereas Remark~\ref{rem:bias} indicates that the magnitude of the bias
grows with $k_n$ as well. Both intuitive properties are visible from
the simulation study provided in the following and from additional
results which we do not show for the sake of brevity.

We begin with a thorough simulation study for a fixed number of
observations $n=22\mbox{,}500$, where we run the simulation $500$ times each.
We have decided to keep the size of the data set fixed in order to work
out the effects that different choices of~$k_n$, or, equivalently, of
$\Delta_n$, have on the finite sample performance of our estimators.
Otherwise, if one fixes $k_n$ or $\Delta_n$ and investigates increasing
sample sizes, it will in general be hard to tell whether a possible
gain in the MSE (say) is due to a more reasonable trade-off between
bias and variance or just to more observations. We briefly discuss
these issues at the end of this section in an additional simulation
study with a fixed number of days $k_n$.

\begin{table}
\tabcolsep=0pt
\caption{Empirical bias and (co)variances of $\sqrt{k_n}
(U_n(\vect{x}) - U(\vect{x}))$ for various choices of $k_n$. Upper five
lines: pure subordinator; lower five lines: subordinator${} + {}$Brownian
motion} \label{tab:U1}
\begin{tabular*}{\textwidth}{@{\extracolsep{\fill}}ld{2.4}cd{2.4}cd{2.4}cccc@{}}
\hline
\multicolumn{1}{@{}l}{$\bolds{x, y}$} &
\multicolumn{2}{c}{$\bolds{2,2}$} & \multicolumn{2}{c}{$\bolds{1,1}$} &
\multicolumn{2}{c}{$\bolds{0.5,0.5}$} & \textbf{2, 0.5} & \textbf{2, 1} & \multicolumn{1}{c@{}}{\textbf{1, 0.5}}
\\[-6pt]
&
\multicolumn{2}{c}{\hrulefill} & \multicolumn{2}{c}{\hrulefill} &
\multicolumn{2}{c}{\hrulefill} & \multicolumn{1}{c}{\hrulefill} &\multicolumn{1}{c}{\hrulefill}& \multicolumn{1}{c@{}}{\hrulefill}\\
\multicolumn{1}{@{}l}{$\bolds{k_n}$} & \multicolumn{1}{c}{$\operatorname{\mathbf{bias}}$} & \multicolumn{1}{c}{$\operatorname{\mathbf{var}}$} &
\multicolumn{1}{c}{$\operatorname{\mathbf{bias}}$} & \multicolumn{1}{c}{$\operatorname{\mathbf{var}}$} &
\multicolumn{1}{c}{$\operatorname{\mathbf{bias}}$} & \multicolumn{1}{c}{$\operatorname{\mathbf{var}}$} & \multicolumn{1}{c}{$\operatorname{\mathbf{cov}}$}&
\multicolumn{1}{c}{$\operatorname{\mathbf{cov}}$} & \multicolumn{1}{c@{}}{$\operatorname{\mathbf{cov}}$} \\
\hline
\phantom{0}50 & -0.0106 & 0.1007 & -0.0077 & 0.1400 & 0.0023 & 0.1915 & 0.0988 &
0.0978 & 0.1376 \\
\phantom{0}75 & -0.0330 & 0.0972 & -0.0229 & 0.1453 & -0.0395 & 0.1956 & 0.1015 &
0.1001 & 0.1435 \\
100 & 0.0168 & 0.1021 & 0.0223 & 0.1375 & 0.0341 & 0.1893 & 0.0996 &
0.0927 & 0.1300 \\
150 & 0.0037 & 0.1061 & 0.0154 & 0.1480 & 0.0470 & 0.2180 & 0.1073 &
0.1106 & 0.1531 \\
250 & 0.0282 & 0.0931 & 0.0547 & 0.1269 & 0.0951 & 0.1845 & 0.0900 &
0.0865 & 0.1245 \\
\phantom{0}50 & -0.0281 & 0.0893 & -0.0120 & 0.1208 & -0.0042 & 0.1863 & 0.0840 &
0.0854 & 0.1233 \\
\phantom{0}75 & 0.0252 & 0.0949 & 0.0115 & 0.1187 & 0.0226 & 0.1861 & 0.0861 &
0.0894 & 0.1216 \\
100 & 0.0126 & 0.0922 &0.0043 & 0.1320 & 0.0401 & 0.1940 & 0.0933 &
0.0932 & 0.1323 \\
150 & -0.0085 & 0.0929 & -0.0127 & 0.1337 & 0.0277 & 0.1991 & 0.0931 &
0.0962 & 0.1371 \\
250 & 0.0128 & 0.1101 & 0.0236 & 0.1395 & 0.0765 & 0.1938 & 0.1049 &
0.1044 & 0.1369 \\
\hline
\end{tabular*}
\end{table}

Despite the fact that we have proven weak convergence of our estimators
in certain function spaces, we restrict ourselves to an analysis of the
finite dimensional properties of our estimators. Let us begin with the
asympotics in Theorem~\ref{prop1} for which we estimate $U(x,x)$ for $x
= 2,1,0.5$. Note that we have $\Cov( \mathbb B(\vect x), \mathbb
B(\vect y)) = (32 \pi)^{-1/2} \approx0.0997$ whenever $\vect x$ or
$\vect y$ equals $(2,2)$, whereas $\Cov( \mathbb B(\vect x), \mathbb
B(\vect y)) = (16 \pi)^{-1/2} \approx0.1410$ if the ``larger'' vector
is $(1,1)$ and finally $\Var( \mathbb B(\vect x)) = (8 \pi)^{-1/2}
\approx0.1995$ for $\vect x = (0.5,0.5)$. Table~\ref{tab:U1} gives
estimated bias and (co)variance for various choices of the number of
trading days, $k_n$. These values are picked in such a way that they
belong to reasonable scenarios in practice. The smallest one, $k_n =
50$, corresponds to $\Delta_n^{-1}=450$ or sampling frequencies of
about a minute, for which microstructure noise already become an issue.
On the other hand, the largest choice of $k_n = 250$ necessitates data
from a process whose jump behavior is homogeneous for quite a long
period of time, namely about one year.

\begin{figure}

\includegraphics{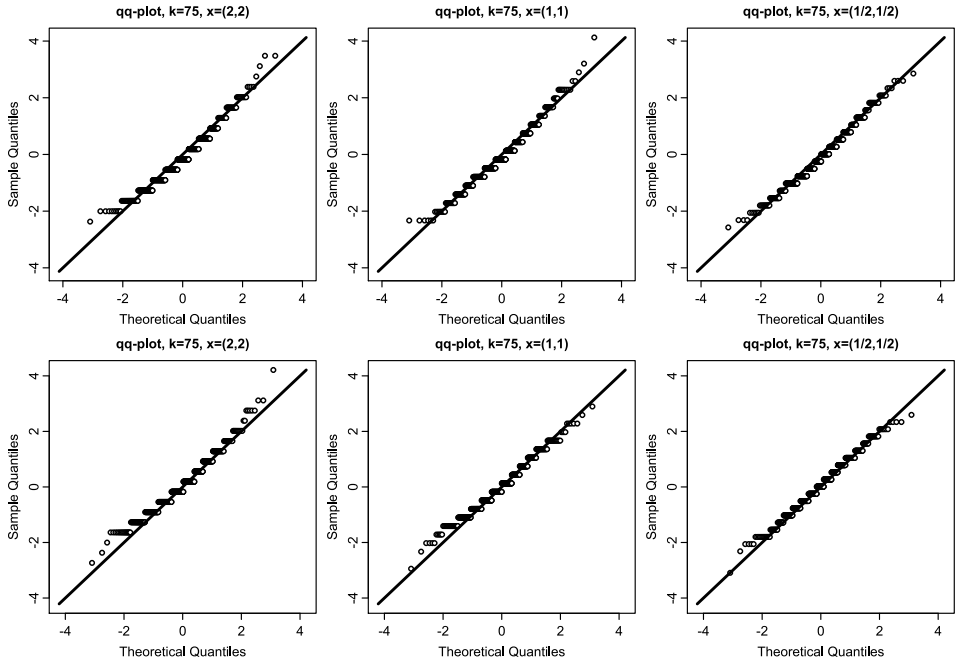}

\caption{QQ-plots of the empirical quantiles of $\sqrt{k_n}
(U_n(\vect
{x}) - U(\vect{x}))$
divided by their sample standard deviation vs. the theoretical
quantiles of the standard normal distribution.
Upper three pictures: pure subordinator; lower three pictures:
subordinator${} + {}$Brownian motion.} \label{fig:qqU}
\end{figure}

Generally, the theoretical (co)variances are well reproduced in both
situations, even though the results look probably a bit better in the
first five lines. This is of course no surprise, since additional
Brownian increments make it harder to infer on the jump measure. In
order to assess how well the normal approximation works apart from bias
and variance, Figure~\ref{fig:qqU} gives QQ-plots for the medium choice
of $k_n = 75$. These plots confirm that the finite sample properties
are indeed satisfying, despite the discrete nature of the test
statistic which simply counts exceedances of certain levels and is
rescaled afterward.

Let us come to the estimation of the Pareto--L\'evy copula. We proceed
in the same way as before and discuss convergence of the finite
dimensional distributions only. For simplicity, we estimate $\Gamma
(x,x)$ for $x = 2,1,0.5$ again, but these are of course different
quantities now. In this case, the variances compute to $\Var(\mathbb
G(\vect x)) = 21/(128x)$, which becomes approximately $0.0820$ for
$x=2$, $0.1641$ for $x=1$, and $0.3281$ for $x=0.5$. Also, for $x>y$ we
have $\Cov( \mathbb G(\vect x), \mathbb G(\vect y)) = 7/32(1/x-\Gamma
(x,y))$. Therefore $\Cov( \mathbb G(\vect2), \mathbb G(\vect{0.5}))
\approx0.0608$, $\Cov( \mathbb G(\vect2), \mathbb G(\vect1))
\approx
0.0718$, and $\Cov( \mathbb G(\vect1), \mathbb G(\vect{0.5}))
\approx
0.1437$. We state their empirical versions in Table~\ref{tab:G1}.

\begin{figure}[b]\vspace*{-3pt}

\includegraphics{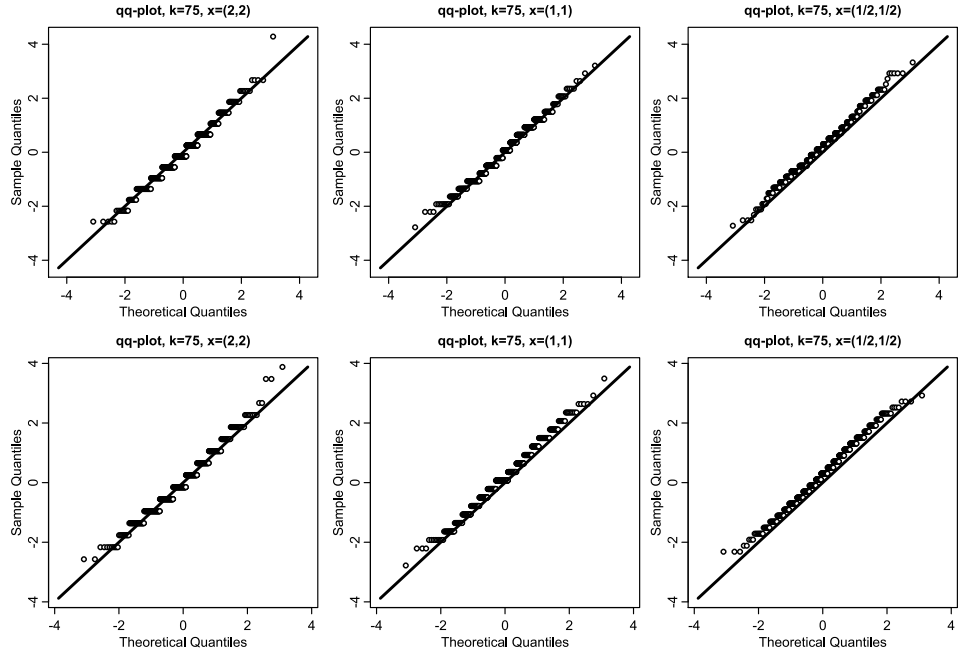}

\caption{QQ-plots of the empirical quantiles of $\sqrt{k_n} (\hat
\Gamma
_n (\vect x) - \Gamma(\vect x))$
divided by their sample standard deviation vs. the theoretical
quantiles of the standard normal distribution.
Upper three pictures: pure subordinator; lower three pictures:
subordinator${} + {}$Brownian motion.} \label{fig:qqG}
\end{figure}

In this case the growth in bias for larger $k_n$ is clearly visible,
and we also have a larger bias when estimating $\Gamma(0.5,0.5)$.
Overall, however, the results are satisfying again, and we see from the
QQ-plot in Figure~\ref{fig:qqG} that the normal approximation works
very well for $k_n = 75$, no matter if a Brownian motion is added or not.

\begin{table}
\tabcolsep=0pt
\def\arraystretch{0.9}
\caption{Empirical bias and (co)variances of $\sqrt{k} (\hat
\Gamma_n (\vect x) - \Gamma(\vect x))$ for various choices of $k_n$.
Upper five lines: pure subordinator; lower five lines: subordinator${}
+{}$Brownian motion} \label{tab:G1}
\begin{tabular*}{\textwidth}{@{\extracolsep{\fill}}ld{2.4}ccccd{1.4}ccc@{}}
\hline
\multicolumn{1}{@{}l}{$\bolds{x, y}$} &
\multicolumn{2}{c}{$\bolds{2,2}$} & \multicolumn{2}{c}{$\bolds{1,1}$} &
\multicolumn{2}{c}{$\bolds{0.5,0.5}$} & \textbf{2, 0.5} & \textbf{2, 1} & \multicolumn{1}{c@{}}{\textbf{1, 0.5}}
\\[-6pt]
&
\multicolumn{2}{c}{\hrulefill} & \multicolumn{2}{c}{\hrulefill} &
\multicolumn{2}{c}{\hrulefill} & \multicolumn{1}{c}{\hrulefill} &\multicolumn{1}{c}{\hrulefill}& \multicolumn{1}{c@{}}{\hrulefill}\\
\multicolumn{1}{@{}l}{$\bolds{k_n}$} & \multicolumn{1}{c}{$\operatorname{\mathbf{bias}}$} & \multicolumn{1}{c}{$\operatorname{\mathbf{var}}$} &
\multicolumn{1}{c}{$\operatorname{\mathbf{bias}}$} & \multicolumn{1}{c}{$\operatorname{\mathbf{var}}$} &
\multicolumn{1}{c}{$\operatorname{\mathbf{bias}}$} & \multicolumn{1}{c}{$\operatorname{\mathbf{var}}$} & \multicolumn{1}{c}{$\operatorname{\mathbf{cov}}$}&
\multicolumn{1}{c}{$\operatorname{\mathbf{cov}}$} & \multicolumn{1}{c@{}}{$\operatorname{\mathbf{cov}}$} \\
\hline
\phantom{0}50 & 0.0141 & 0.0827 & 0.0455 & 0.1740 & 0.0863 & 0.3520 & 0.0777 &
0.0668 & 0.1599 \\
\phantom{0}75 & -0.0082 & 0.0874 & 0.0173 & 0.1653 & 0.1252 & 0.3428 & 0.0740 &
0.0690 & 0.1459 \\
100 & 0.0502 & 0.0783 & 0.0894 & 0.1708 & 0.1748 & 0.3400 & 0.0685 &
0.0508 & 0.1547 \\
150 & 0.0356 & 0.0862 & 0.1182 & 0.1646 & 0.3176 & 0.3324 & 0.0744 &
0.0698 & 0.1421 \\
250 & 0.0505 & 0.0732 & 0.1493 & 0.1605 & 0.5625 & 0.313 & 0.0678 &
0.0568 & 0.1390 \\
\phantom{0}50 & 0.0345 & 0.0790 & 0.0560 & 0.1637 & 0.1021 & 0.3263 & 0.0699 &
0.0639 & 0.1389 \\
\phantom{0}75 & 0.0091 & 0.0886 & 0.0753 & 0.1760 & 0.1522 & 0.3508 & 0.0832 &
0.0729 & 0.1522 \\
100 & 0.0312 & 0.0745 & 0.0776 & 0.1530 & 0.1480 & 0.3033 & 0.0610 &
0.0558 & 0.1305 \\
150 & 0.0284 & 0.0866 & 0.0988 & 0.1694 & 0.2074 & 0.3337 & 0.0746 &
0.0725 & 0.1486 \\
250 & 0.0602 & 0.0762 & 0.1515 & 0.1528 & 0.3012 & 0.3452 & 0.0675 &
0.0545 & 0.1455 \\
\hline
\end{tabular*}   \vspace*{-3pt}
\end{table}

Finally, let us briefly discuss the performance of our estimators in
case of a fixed number $k_n$ and for increasing sampling frequencies
$\Delta_n^{-1}$. We restrict ourselves to the case of estimating $U$
for a pure subordinator, as other settings lead to similar results. For
$k_n$ fixed, one would expect that choosing a rather low frequency
yields the worst results since in this case the bias is the largest,
or, equivalently, the jumps are most difficult to detect. On the other
hand, for growing $\Delta_n^{-1}$ this bias becomes smaller, but
otherwise not much is to be gained, as the sampling frequency does not
affect the rate of convergence. The results in Table~\ref{tab:01},
where we consider all possible combinations of $k_n, \Delta_n^{-1}\in
\{
50,100,150,200\}$, support these findings. For financial applications
this means that, provided sufficient data is available, a~trade-off has
to be made: on the one hand, one should choose the frequency as high as
possible; but on the other hand, due to microstructure noise issues,
care is needed regarding a possibly oversized frequency.


\begin{table}
\tabcolsep=0pt
\caption{Empirical bias and (co)variances of $\sqrt{k_n}
(U_n(\vect{x}) - U(\vect{x}))$ for various choices of $k_n$ and~$\Delta
_n^{-1}$. In all cases: pure subordinator} \label{tab:01}
\begin{tabular*}{\textwidth}{@{\extracolsep{\fill}}ld{2.4}cd{2.4}cd{2.4}cccc@{}}
\hline
\multicolumn{1}{@{}l}{$\bolds{x, y}$} &
\multicolumn{2}{c}{$\bolds{2,2}$} & \multicolumn{2}{c}{$\bolds{1,1}$} &
\multicolumn{2}{c}{$\bolds{0.5,0.5}$} & \textbf{2, 0.5} & \textbf{2, 1} & \textbf{1, 0.5}
\\[-6pt]
&
\multicolumn{2}{c}{\hrulefill} & \multicolumn{2}{c}{\hrulefill} &
\multicolumn{2}{c}{\hrulefill} &\multicolumn{1}{c}{\hrulefill}  &\multicolumn{1}{c}{\hrulefill}& \multicolumn{1}{c@{}}{\hrulefill}\\
\multicolumn{1}{@{}l}{$\bolds{\Delta_n^{-1}}$} & \multicolumn{1}{c}{$\operatorname{\mathbf{bias}}$} & \multicolumn{1}{c}{$\operatorname{\mathbf{var}}$} &
\multicolumn{1}{c}{$\operatorname{\mathbf{bias}}$} & \multicolumn{1}{c}{$\operatorname{\mathbf{var}}$} &
\multicolumn{1}{c}{$\operatorname{\mathbf{bias}}$} & \multicolumn{1}{c}{$\operatorname{\mathbf{var}}$} & \multicolumn{1}{c}{$\operatorname{\mathbf{cov}}$}&
\multicolumn{1}{c}{$\operatorname{\mathbf{cov}}$} & \multicolumn{1}{c@{}}{$\operatorname{\mathbf{cov}}$} \\
\hline
\multicolumn{10}{c}{$k_n=50$} \\
\phantom{0}50 & 0.0338 & 0.1065 & 0.0285 & 0.1425 & 0.0648 & 0.2072 & 0.1037 &
0.1017 & 0.1407 \\
100 & 0.0010 & 0.0976 & 0.0005 & 0.1316 & 0.0354 & 0.1972 & 0.0942 &
0.0936 & 0.1336 \\
150 & 0.0007 & 0.1058 & 0.0138 & 0.1445 & 0.0204 & 0.2038 & 0.1044 &
0.1043 & 0.1446 \\
200 & 0.0027 & 0.1005 & 0.0059 & 0.1370 & 0.0286 & 0.1928 & 0.0991 &
0.0935 & 0.1343 \\[3pt]
\multicolumn{10}{c}{$k_n=100$} \\
\phantom{0}50 & 0.0194 & 0.1032 & 0.0589 & 0.1424 & 0.1071 & 0.2052 & 0.1024 &
0.1026 & 0.1421 \\
100 & -0.0200 & 0.1007 & 0.0017 & 0.1466 & 0.0175 & 0.1981 & 0.1034 &
0.1025 & 0.1434 \\
150 & 0.0150 & 0.1002 & 0.0175 & 0.1480 & 0.0275 & 0.1999 & 0.1021 &
0.0979 & 0.1444 \\
200 & 0.0202 & 0.0957 & 0.0103 & 0.1326 & 0.0299 & 0.2059 & 0.0953 &
0.1017 & 0.1392 \\[3pt]
\multicolumn{10}{c}{$k_n=150$} \\
\phantom{0}50 & 0.0312 & 0.1034 & 0.0659 & 0.1450 & 0.1286 & 0.2080 & 0.1037 &
0.1057 & 0.1437 \\
100 & 0.0132 & 0.0988 & 0.0205 & 0.1369 & 0.0489 & 0.2092 & 0.0960 &
0.0984 & 0.1407 \\
150 & 0.0037 & 0.1061 & 0.0154 & 0.1480 & 0.0470 & 0.2180 & 0.1073 &
0.1106 & 0.1531 \\
200 & -0.0098 & 0.0948 & -0.0122 & 0.1315 & -0.0156 & 0.1891 & 0.0923
& 0.0884 & 0.1299 \\[3pt]
\multicolumn{10}{c}{$k_n=200$} \\
\phantom{0}50 & 0.0134 & 0.0990 & 0.0468 & 0.1334 & 0.1237 & 0.1947 & 0.0932 &
0.0892 & 0.1306 \\
100 & -0.0012 & 0.0995 & 0.0136 & 0.1415 & 0.0554 & 0.2092 & 0.0986 &
0.1003 & 0.1437 \\
150 & -0.0006 & 0.0966 & 0.0217 & 0.1371 & 0.0413 & 0.1966 & 0.0944 &
0.0950 & 0.1372 \\
200 & 0.0011 & 0.0934 & 0.0292 & 0.1371 & 0.0305 & 0.2074 & 0.0942 &
0.0961 & 0.1416 \\
\hline
\end{tabular*}
\end{table}

\subsection{Illustration}
In the present section we are going to apply the estimation techniques
developed in the previous sections to infer on the jump dependence of
some specific financial data set.

More precisely, we consider the logarithm of one-minute Nasdaq stock
prices of Apple Inc. and Microsoft Corporation in the third quarter of
2012, which consists of $62$ trading days. After some cleaning of the
corresponding time series we obtain a two-dimensional data sample of
increments (log returns) of size $n=21\mbox{,}864$. From the simulation
results in Section~\ref{subsec:sim} we know that the choice $k_n=62$
yields a reasonable trade-off between bias and variance.

Due to missing observations or errors in data, we formally apply the
procedure from Section~\ref{sec:irr2}, though the observations are in
principle quite close to a regular sampling scheme. We have chosen a
frequency of one minute returns in order to have a large amount of data
while not being too much affected by microstructure effects which our
method does not correct for; see Section~\ref{Noise}. Note that the
results in this paper are stated for observations of a bivariate
semimartingale with a constant L\'evy measure which is probably a too
simple model for a bivariate price process. What might be less
restrictive, is to assume a time-homogeneous Pareto--L\'evy copula, if
one is only interested in the dependence structure of the assets.

\begin{figure}[b]

\includegraphics{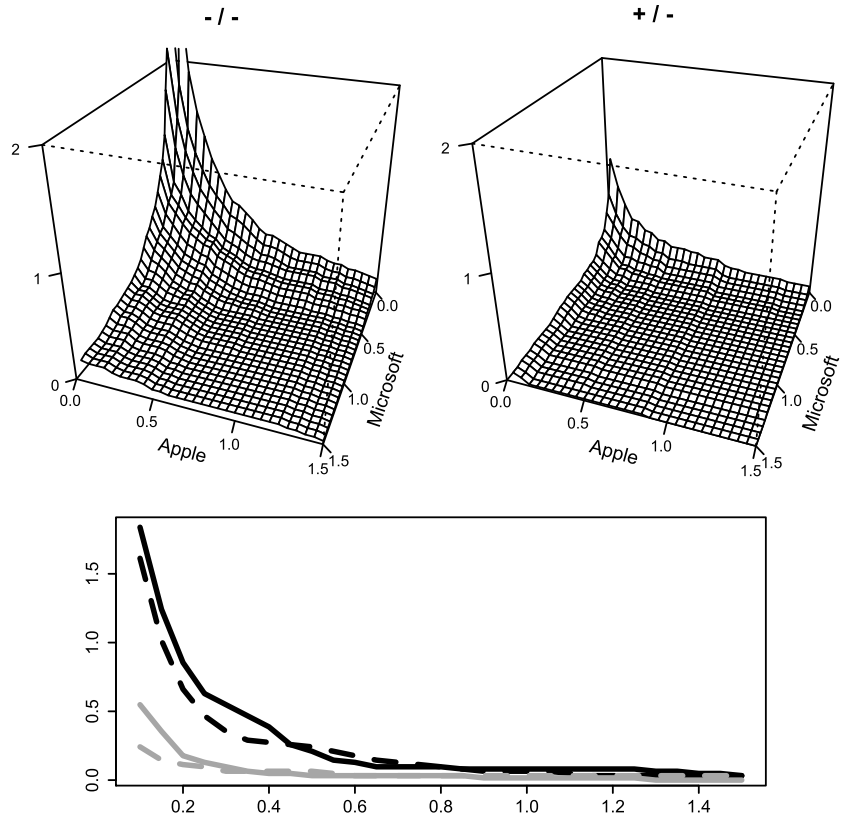}

\caption{Upper two pictures: empirical Pareto--L\'{e}vy copulas for the
$(-/-)$ and the $(+/-)$ dependence.
Lower picture: empirical Pareto--L\'{e}vy copulas along the diagonal
for $(-/-)$ (black solid line),
$(+/+)$ (black dashed), $(+/-)$ (gray solid) and $(-/+)$ (gray dashed)
dependence.}\label{pic:appmic}
\end{figure}

For this reason, we are interested in margin-free estimates of the jump
dependence only, which means that we restrict ourselves to the
estimation of the Pareto--L\'evy copula $\Gamma$. There are four
possible types of jump dependence: positive jumps in both components
($+/+$), positive jumps in the log-returns of Apple may be associated
with negative jumps in the log-returns of Microsoft ($+/-$), and vice
versa ($-/+$), and negative jumps in both components ($-/-$). To obtain
estimates in the latter three cases, we simply use\vspace*{1pt} negative log-returns
in the corresponding components in the definition of~$\hat\Gamma_n$.

The results are depicted in Figure~\ref{pic:appmic}. In the two upper
pictures, we plot the graph of the empirical Pareto--L\'evy copula on
the set $[0.05, 1.5]^2$ for the dependencies $(-/-)$ and $(+/-)$,
respectively. The corresponding graphs for $(+/+)$ and $(-/+)$,
respectively, look very similar and are therefore omitted. In the lower
picture, we plot the restriction of the graphs of the empirical
Pareto--L\'evy copula to the main diagonal for all four kinds of dependence.

The following are the main findings:
\begin{itemize}
\item The dependence between positive or negative jumps in both
components ($+/+$ and $-/-$) is generally much stronger than the
dependence between positive and negative jumps ($+/-$ and $-/+$).
\item Comparing the $(+/+)$ and the $(-/-)$ dependence, the latter is
slightly stronger.
\item Comparing the $(+/-)$ and the $(-/+)$ dependence, we observe that
it is more likely that positive jumps in the Apple returns occur
simultaneously with negative jumps in the Microsoft returns than vice
versa. Both are, however, close to being independent; that is, the
estimated Pareto--L\'evy copula is close to $\Gamma_\perp$ from
Proposition~\ref{prop:plc}.
\end{itemize}

\section{Conclusions} \label{sec:conc}

In this paper we have investigated the problem of estimating both the
bivariate L\'evy measure and the (Pareto) L\'evy copula in a
nonparametric way. Our estimators are based on counting joint large
increments of a bivariate L\'evy process, and in both cases we were
able to prove weak convergence in appropriate function spaces. An
extension to the case of irregular and/or asynchronous observations is
provided as well. What still remains an open problem is the development
of similar methods when microstructure noise is present, as indicated
in Section~\ref{Noise}.

From a statistical point of view, it might also be interesting to
construct several nonparametric tests concerning the dependence
structure of a multivariate It\^o semimartingale. Using the methods
from this work, one should be able to check first whether the entire
jump measure (or just the jump dependence) is indeed constant over
time, while under the assumption of a genuine L\'evy jump part these
procedures could include estimation of certain functionals of $\Gamma$
or $U$, as well as tests for independence or tests for a parametric
form of these functions. For this reason, it would be important to
establish a thorough theory concerning (Pareto) L\'evy copulas which
relates functionals of $\Gamma$ to certain dependence properties, as in
the case of ordinary copulas for which standard measures such as
Kendall's $\tau$ or Spearman's $\rho$ can be written as integrals over
$C$ and are thus accessible through nonparametric estimation of the copula.

\begin{appendix}
\section*{Appendix} \label{sec:aux}

In this section we present the proofs of the main Theorems~\ref{prop1}
and~\ref{theo1}. Proofs of the additional results as well as some
technical lemmas are postponed to a supplementary Appendix.

\subsection{\texorpdfstring{Proof of Theorem \protect\ref{prop1}}{Proof of Theorem 4.2}}
Before we begin with
the proof, note that, due to Theorem 1.6.1 in \citet{vandwell1996}, weak
convergence in $\mathcal B_\infty((0, \infty]^2)$ is equivalent to weak
convergence on each $\ell^{\infty}(T_k)$, which is the space of all
bounded functions on $T_k$ endowed with the uniform norm. Therefore, it
is possible to fix one such $T_k$ throughout the rest of the proof.

Let us introduce some additional notation. We define a class of
functions $\mathcal F_n = \{f_{n,\vect{x}}\dvtx \vect{x} \in T_k\}$ via
\[
f_{n,\vect{x}}(\vect{p}) = \sqrt{{n}/{k_n}} ( 1_{\{ \vect{p}
\ge
\vect{x} \ge(0,0) \} } +
1_{ \{ p_1 \ge x_1, x_2=-\infty\} } + 1_{\{
p_1 \ge x_2, x_1=-\infty\} } ).
\]
Furthermore, we set
\[
\Bga_n(\vect{x})=\sqrt{k_n} \bigl( U_n(
\vect{x}) - \E\bigl[U_n(\vect{x})\bigr] \bigr) = \frac{1}{\sqrt n}
\sum_{j=1}^n \bigl(f_{n,\vect
{x}}\bigl(
\Delta _j^n \vect X\bigr) - \E\bigl[f_{n,\vect{x}}
\bigl(\Delta_j^n \vect X\bigr)\bigr] \bigr).
\]
As a consequence of Lemma~\ref{lem1}, it is sufficient to discuss weak
convergence of $\Bga_n(\vect{x})$ only. Indeed, let $\vect{x} \in T_k$.
Then by stationarity of increments of $\vect X$ and using $k_n = n
\Delta_n$, we have
\begin{eqnarray*}
&&\E\bigl[U_n(\vect{x})\bigr] - U(\vect{x})
\\
&&\qquad= \Delta_n^{-1} \PP \bigl( \Delta _1^n
X^{(1)} \geq x_1, \Delta_1^n
X^{(2)} \geq x_2 \bigr) - \nu \bigl([x_1,\infty
) \times[x_2,\infty)\bigr).
\end{eqnarray*}
This quantity is bounded by $K\Delta_n$ due to Lemma~\ref{lem1}, so the
growth condition $\sqrt{k_n} \Delta_n \to0$ ensures that $\sqrt{k_n}
(\gamma_n(\vect{x}) - \Bga_n(\vect{x}))$ is uniformly small on each
fixed~$T_k$.

In order to prove $\Bga_n(\vect{x}) \weak{\mathbb B}(\vect{x})$ on
$\ell^{\infty}(T_k)$ we will employ Theorem 11.20 in \citet{kosorok2008}
for which several intermediate results have to be shown. To begin with, set
$
F_{n}(\vect{p}) = \sqrt{{n}/{k_n}} 1_{\{\vect{p} \in T_k\}},
$
which is a sequence of integrable (with respect to any probability
measure) envelopes. The first two steps are related to the class of
functions $\mathcal F_n$. We start with the proof of an entropy
condition, namely
\[
\limsup_{n \to\infty} \sup_Q \int
_0^1 \sqrt{\log N\bigl(\eps\| F_n \|
_{Q,2}, \mathcal F_n, L_2(Q)\bigr)} \,d\eps<
\infty, \label{entr}
\]
where $N$ denotes the covering number of the set $\mathcal F_n$, and
the supremum runs over all probability measures $Q$ with finite support
such that $\| F_n \|_{Q,2} =  (\int F^2_n(\vect{p}) \,dQ(\vect
{p}) )^{1/2} > 0$. Thanks to the special form of $\mathcal F_n$,
this result is a simple consequence of Lemma 11.21 in \citet
{kosorok2008}: it suffices to check that each $\mathcal F_n$ is a
VC-class with VC-index $5$.
This follows from the fact that each finite subset of $\HH$ of size 5
has either a subset of 3 elements in $[0,\infty]^2\setminus\{(0,0)\}$,
or a subset of two elements in one of the stripes through $-\infty$. In
neither of the cases these subsets can be shattered by the sets deduced
from the indicators in the definition of $f_{n,\vect{x}}$.

The second condition to check is that $\mathcal F_n$ is almost
measurable Suslin, and it follows from Lemma 11.15\vadjust{\goodbreak} and the discussion
on page 224 in \citet{kosorok2008} that it is sufficient to prove
separability of $\mathcal F_n$, that is, to show the existence of a
countable subset $T_{n,k}$ of $T_k$ such that
\[
\mathbb P^* \Bigl(\sup_{\vect{x} \in T_k} \inf_{\vect{y} \in T_{n,k}}
\bigl\llvert f_{n,\vect{x}}\bigl(\Delta_j^n \vect X
\bigr) - f_{n,\vect{y}
}\bigl(\Delta _j^n \vect X\bigr)
\bigr\rrvert > 0 \Bigr) = 0.
\]
Here, $\mathbb P^*$ denotes the outer expectation, since measurability
of the event within the brackets is not ensured. Set $T_{n,k} = T_k
\cap\overline{\mathbb Q}^2$. Then, for each $\omega$ and each $\vect
{x} \in T_k$, there exists a $\vect{y} \in T_{n,k}$ such that
$f_{n,\vect{x}}(\Delta_j^n \vect X(\om)) = f_{n,\vect{y}}(\Delta_j^n
\vect X(\om))$, since the $f_{n,\vect{x}}$ are indicator functions.
This proves separability of $\mathcal F_n$.

The remaining steps regard the behavior of the variances and
covariances of the $f_{n,\vect{x}}$ and their envelopes. We have
%
%
\begin{equation}
\label{cov} \lim_{n \to\infty} \E\bigl[\Bga_n(\vect{x})
\Bga_n(\vect{y})\bigr] = \lim_{n \to
\infty} \E
\bigl[f_{n,\vect{x}}\bigl(\Delta_j^n \vect X\bigr)
f_{n,\vect
{y}}\bigl(\Delta_j^n \vect X\bigr)\bigr] =
U(\vect{x} \vee\vect{y})
\end{equation}
as well as
\[
\lim_{n \to\infty} \E\bigl[F_{n}^2\bigl(
\Delta_j^n \vect X\bigr)\bigr] \le U(1/k,-\infty) + U(-
\infty,1/k)
\]
and
\[
\lim_{n \to\infty} \E\bigl[F_{n}^2\bigl(
\Delta_j^n \vect X\bigr) 1_{\{
F_{n}(\Delta
_j^n \vect X) > \eps\sqrt n \}}\bigr] \leq\lim
_{n \to\infty} \E \bigl[F_{n}^2\bigl(
\Delta_j^n \vect X\bigr)\bigr] (\eps\sqrt{k_n}
)^{-1} \to0.
\]
Finally, as in (\ref{cov}) we have for $\vect{x}, \vect{y} \in T_k$ that
\begin{eqnarray*}
\rho_n(\vect x, \vect y) &=& \bigl(\E \bigl[ \bigl(f_{n,\vect
x}
\bigl(\Delta _j^n \vect X\bigr) - f_{n,\vect y}\bigl(
\Delta_j^n \vect X\bigr) \bigr)^2 \bigr]
\bigr)^{1/2}
\\
&\to& \bigl(U(\vect x)+ U(\vect y) - 2U(\vect x \vee\vect y)
\bigr)^{1/2} = \rho(\vect x, \vect y),
\end{eqnarray*}
and due to Lemma~\ref{lem1} the convergence holds uniformly as well.
This completes the proof. 

\subsection{\texorpdfstring{Proof of Theorem \protect\ref{theo1}}{Proof of Theorem 4.6}}
Let $\mathcal B_\infty^0((0,\infty]^2)\subset\mathcal B_\infty
((0,\infty
]^2)$ and\break $\mathcal B_\infty^0((0,\infty]) \subset\mathcal B_\infty
((0,\infty])$ denote the space of all tail integrals of bivariate L\'
evy measures concentrated on the first quadrant or of univariate L\'evy
measures concentrated on $(0,\infty]$, respectively.
Consider the mapping $\Phi\dvtx\break  \mathcal B^0_\infty((0,\infty]^2) \times
(\mathcal B^0_\infty((0,\infty]))^2 \to\mathcal B_\infty((0,\infty
]^2)$, defined by $\Phi= \Phi_3 \circ\Phi_2 \circ\Phi_1$ with
\begin{eqnarray*}
&&\Phi_1\dvtx\quad  \mathcal B^0_\infty\bigl((0,
\infty]^2\bigr) \times\bigl(\mathcal B^0_\infty
\bigl((0,\infty]\bigr)\bigr)^2 \to\mathcal B^0_\infty\bigl(
(0,\infty]^2\bigr) \times\bigl(\mathcal B^-_\infty \bigl((0,
\infty]\bigr)\bigr)^2
\\
&&\phantom{\Phi_1\dvtx\quad}(U, U_1,U_2) \mapsto\bigl(U, U_1^-,
U_2^-\bigr),
\\
&&\Phi_2\dvtx\quad  \mathcal B^0_\infty\bigl((0,
\infty]^2\bigr) \times\bigl(\mathcal B^-_\infty \bigl((0,\infty]\bigr)
\bigr)^2 \to\mathcal B^0_\infty\bigl((0,
\infty]^2\bigr) \times\bigl(\mathcal B^p_\infty
\bigl([0,\infty)\bigr)\bigr)^2
\\
&&\phantom{\Phi_1\dvtx\quad}(U, V_1,V_2) \mapsto(U, V_1\circ P,
V_2 \circ P),
\\
&&\Phi_3\dvtx\quad \mathcal B^0_\infty\bigl((0,
\infty]^2\bigr) \times\bigl(\mathcal B^p_\infty
\bigl([0,\infty)\bigr)\bigr)^2 \to\mathcal
B_\infty\bigl((0,\infty]^2\bigr)
\\
&&\phantom{\Phi_1\dvtx\quad}(U,G_1,G_2) \mapsto U({G_1},
{G_2}), 
\end{eqnarray*}
where $P(x)=1/x$ and where, in the last step, $G_i(\infty)=\infty$. %
Moreover,\break $\mathcal B^-_\infty((0,\infty])\subset\mathcal B_\infty
((0,\infty])$ and $\mathcal B^p_\infty([0,\infty))\subset\mathcal
B_\infty([0,\infty))$ are\vadjust{\goodbreak} defined as the images of the associated
function spaces under the respective mappings. Set also $\widetilde
\Gamma_{n,1}(x)=U_n(U_1^-(1/x), -\infty)$ and $\widetilde\Gamma
_{n,2}(x)=U_n(-\infty, U_2^-(1/x))$. The proof will now basically
consist of two steps. We start with discussing weak convergence of
%
%
\begin{equation}
\label{eq:hadam} \sqrt{k_n} \bigl( \Phi(\widetilde
\Gamma_n, \widetilde\Gamma_{n,1}, \widetilde
\Gamma_{n,2}) - \Phi(\Gamma, P,P) \bigr) \weak\mathbb G,
\end{equation}
whereas this result is transferred to the original claim later on.

Let us begin with the proof of (\ref{eq:hadam}). This assertion follows
from the functional delta method in topological vector spaces [see
\citet
{vandwell1996}], if we prove first that
\[
\sqrt{k_n} \bigl\{ (\widetilde\Gamma_n, \widetilde
\Gamma_{n,1}, \widetilde\Gamma_{n,2}) - (\Gamma,P,P) \bigr\}
\weak\bigl(\widetilde{\mathbb G}, \widetilde{\mathbb G}(\cdot,-\infty),
\widetilde{\mathbb G}(-\infty,\cdot)\bigr)
\]
in $\mathcal B_\infty((0,\infty]^2) \times(\mathcal B_\infty
((0,\infty
]))^2$ and second that $\Phi$ is Hadamard-differentiable at $(\Gamma,P,P)$ tangentially to suitable subspaces with derivative
%
%
\begin{equation}
\label{eq:deriv}\qquad \bigl( \Phi'_{(\Gamma,P,P)}(U, U_1,
U_2) \bigr) (\vect{u}) = U(\vect u) + u_1^2
\dot\Gamma_1(\vect u) U_1(u_1) +
u_2^2 \dot \Gamma_2(\vect u)
U_2(u_2),
\end{equation}
where the summands involving the partial derivatives on the right-hand
side are defined as $0$ if one of the coordinates of $\vect u$ equals
$\infty$.
The first claim follows easily from Theorem~\ref{prop1} and the
continuous mapping theorem. Regarding the second assertion we need to
clarify the metrics on the corresponding spaces. The canonical
definitions are
\[
d(f,g) = \sum_{k=1}^\infty2^{-k} \bigl(
\| f-g \|_{T_k} \wedge1 \bigr),
\]
where $T_k=[1/k,\infty]^2$ in case of $\mathcal B_\infty((0,\infty
]^2)$, while $T_k=[1/k,\infty]$ and $T_k=[0,k]$ for $\mathcal B_\infty
((0,\infty])$ and $\mathcal B_\infty([0,\infty))$, respectively.
Unfortunately, the mapping $\Phi_1$ is not Hadamard-differentiable with
respect to these metrics [see the proof of Lemma A.2 in
the supplementary material \citet{buecvett2013}], whence we need to
consider the weaker modifications
\[
d_2(f,g) = \sum_{k=1}^\infty2^{-k}
\bigl( \| f-g \|_{S_k} \wedge1 \bigr),
\]
where $S_k=([1/k,k]\cup\{\infty\} )^2$ in case of $\mathcal B_\infty
((0,\infty]^2)$, while $S_k=[1/k,k]\cup\{\infty\}$ and $S_k=\{0\}
\cup
[1/k,k]$ for $\mathcal B_\infty((0,\infty])$ and $\mathcal B_\infty
([0,\infty))$, respectively. With these modifications, it follows from
Lemma A.1 in the supplementary material \citet{buecvett2013}
and the chain rule that
\[
\Phi\dvtx\quad \bigl(\mathcal B^0_\infty\bigl((0,\infty]^2\bigr),d
\bigr) \times\bigl(\mathcal B^0_\infty \bigl((0,\infty]\bigr),d
\bigr)^2 \to\bigl(\mathcal B^0_\infty\bigl((0,
\infty]^2\bigr),d_2\bigr)
\]
is Hadamard-differentiable at $(\Gamma,P,P)$ with derivative as
specified in \eqref{eq:deriv} tangentially to
%
%
\begin{eqnarray}
\label{eq:tang} \mathbb D_{0} &=& \Bigl\{ (U,U_1,U_2)
\in\mathcal C\bigl((0,\infty]^2\bigr) \times \bigl(\mathcal C\bigl((0,
\infty]\bigr)\bigr)^2 \mid U_j(\infty)=0,
\nonumber
\\[-8pt]
\\[-8pt]
\nonumber
&&\hspace*{173pt} \lim
_{x\to0} x^2U_j(x)=0 \Bigr\}.
\end{eqnarray}
Here, $\mathcal C((0,\infty]^2)$ and $C((0,\infty])$ denote the set of
all functions on $(0,\infty]^2$ and $(0,\infty]$ that are continuous
with respect to the pseudo metrics $\rho(\vect u, \vect v)=|\Gamma
(\vect u)-\Gamma(\vect v)|^{1/2}$ and $\rho(u,v)=|1/u-1/v|^{1/2}$,
respectively.\vspace*{1pt} Hence, observing $(\widetilde{\mathbb G}, \widetilde
{\mathbb G}(\cdot,-\infty), \widetilde{\mathbb G}(-\infty,\cdot
))\in
\mathbb D_0$, the functional delta method yields
\[
\sqrt{k_n} \bigl( \Phi(\widetilde\Gamma_n, \widetilde
\Gamma_{n,1}, \widetilde\Gamma_{n,2}) - \Phi(\Gamma, P,P) \bigr)
\weak\mathbb G
\]
%
in $(\mathcal B_\infty((0,\infty]^2),d_2)$.

We will use the approximation Theorem 4.2 in \citet{billingsley1968},
adapted to the concept of weak convergence in the sense of
Hoffmann--J\o rgensen, to transfer this result to weak convergence in
$(\ell^\infty([\eta, \infty]^2), \|\cdot\|_\infty)$ for all $\eta>0$
and hence in $(\mathcal B_\infty((0,\infty]^2),d)$. To this end, define
\[
W_n(\vect u) =\sqrt{k_n}\bigl(\Phi(\widetilde
\Gamma_n, \widetilde\Gamma _{n,1}, \widetilde
\Gamma_{n,2}) - \Phi(\Gamma, P,P)\bigr) (\vect u)\vspace*{-1pt}
\]
and
\[
W_{n,M}(\vect u) = \sqrt{k_n}\bigl(\Phi(\widetilde
\Gamma_n, \widetilde \Gamma_{n,1}, \widetilde
\Gamma_{n,2}) - \Phi(\Gamma, P,P)\bigr) (\vect u) 1_{\{ \vect u \in[\eta,M]^2 \}}.
\]
Then $W_{n,M}(\vect u) \weak\mathbb G_M (\vect u):= \mathbb G(\vect
u) 1_{\{ \vect u \in[\eta,M]^2 \}}$ for $n\to\infty$ and $\mathbb G_M
(\vect u) \weak\mathbb G(\vect u)$ for $M\to\infty$ in $(\ell
^\infty
([\eta, \infty])^2, \|\cdot\|_\infty)$, and it remains to prove that
\[
\limsup_{n\to\infty} \mathbb P^* \Bigl( \sup_{u_1>M \ \mathrm{or}\ u_2>M}
\bigl\llvert \sqrt{k_n}\bigl(\Phi(\widetilde\Gamma_n,
\widetilde\Gamma_{n,1}, \widetilde\Gamma_{n,2}) - \Phi(\Gamma,
P,P)\bigr) (\vect u) \bigr\rrvert >\eps \Bigr)
\]
converges to $0$ for $M\to\infty$.
Noting that $\Phi(\Gamma, P,P) = \Gamma$, the probability can be
bounded by
\begin{eqnarray*}
&&\mathbb P^* \Bigl( \sup_{u_1 \ge M/2\ \mathrm{or}\ u_2 \ge M/2} \bigl| \sqrt
{k_n}(\widetilde\Gamma_n - \Gamma) (\vect u) \bigr| >\eps
\Bigr)
\\
&&\qquad{}+ \mathbb P^* \bigl( \exists\vect u \mbox{ with } u_1> M \mbox{ or }
u_2> M\dvtx \widetilde\Gamma_{ni}^-\circ
P(u_i) < M/2, i=1,2 \bigr).
\end{eqnarray*}
The Portmanteau theorem implies that the $\limsup$ of the first
probability converges to $0$ for $M \to\infty$ using Theorem \ref
{prop1}. Furthermore, some thoughts reveal that $\widetilde\Gamma
_{n,i}^-(z) = 1/(U_i(U^-_{n,i}(z)))$ for all $z > 0$. Due to
monotonicity of $\widetilde\Gamma_{ni}^-\circ P$, the second
probability is bounded by $\mathbb P( \widetilde\Gamma_{ni}^-\circ
P(M)\le M/2, i=1,2)$, which thus converges to $0$ for $n\to\infty$
observing that $\widetilde\Gamma_{ni}^-\circ P(M)=M+o_P(1)$.

In the final step we will prove $\sqrt{k_n} ( \hat\Gamma_n - \Gamma)
\weak\mathbb G$ in each $(\ell^\infty([\eta, \infty]^2),\break  \|\cdot\|
_\infty)$, for which we heavily rely on the fact that the same result
holds for the statistic discussed above. As a consequence of the
identity $\widetilde\Gamma_{n,i}^-(z) = 1/(U_i(U^-_{n,i}(z)))$, we
have that $\Phi(\widetilde\Gamma_n, \widetilde\Gamma_{n,1},
\widetilde\Gamma_{n,2})(\vect u)$ and $\hat\Gamma_n(\vect u)$
coincide as long as $U^-_{n,i}(1/u_i) \neq0$ for $i=1,2$. By
monotonicity, it is therefore sufficient to prove that the probability
of $U^-_{n,i}(1/\eta) = 0$ becomes small, which is precisely
\[
\lim_{n \to\infty} \PP \bigl( U^-_{n,i}(1/\eta) = 0 \bigr) =
0.
\]
To this end, let $N_i(n)$ denote the number of positive increments of
$X^{(i)}$. By definition of the generalized inverse function in (\ref
{geni}) we have that $U^-_{n,i}(1/\eta) = 0$ is equivalent to $1/\eta
\ge N_i(n)/k_n$ or $N_i(n) \leq k_n/\eta$. Furthermore, letting
$M_i(n)$ be the number of positive increments of the process $Z^{(i)}_t
= a_i t + B^{(i)}_t$, we see that it is sufficient to prove
\[
\lim_{n \to\infty} \PP \bigl( M_i(n) \leq
k_n/\eta \bigr) = 0,
\]
since $\vect X$ does not admit negative jumps. Note that we have
\[
\PP \bigl(\Delta_j^n Z^{(i)} > 0 \bigr) = \PP
\bigl(\Delta_j^n B^{(i)} > -a_i
\Delta_n \bigr) = \PP \bigl(N > -a_i \Delta^{1/2}_n
\bigr) ={1}/{2} + o(1),
\]
where $N$ is a standard Gaussian variable. Let $n$ be large enough in
order for the probability above to be larger than 1/3. For such $n$, we
conclude easily that
\[
\PP \bigl( M_i(n) \leq k_n/\eta \bigr) \leq\PP \bigl(
\mbox {Bin}(n,1/3) \leq k_n/\eta \bigr) \to0,
\]
for example, from Markov inequality and (\ref{eq:cond}). This completes the
proof.
\end{appendix}

\section*{Acknowledgments}
The authors would like to thank two unknown
referees and an Associate Editor for their constructive comments on an
earlier version of this manuscript, which led to a substantial
improvement of the paper.

\begin{supplement}[id=suppA]
\stitle{Proof of auxiliary results}
\slink[doi]{10.1214/13-AOS1116SUPP} 
\sdatatype{.pdf}
\sfilename{aos1116\_supp.pdf}
\sdescription{In this supplement we present the proofs of the
remaining results from the main corpus as well as two lemmas which are
used in the proof of Theorem~\ref{theo1}.}
\end{supplement}

\printaddresses

\end{document}